\documentclass[a4paper,11pt]{article}                                                   

\usepackage[english]{babel}
\usepackage{amsfonts,amsmath,amsthm,graphicx,a4wide}   
\usepackage[colorlinks=true]{hyperref}

\usepackage{hyperref}

\usepackage{graphicx} 
\usepackage{amsmath} 
\usepackage{amssymb}  
\usepackage{mathrsfs}
\usepackage{bbm}
\usepackage{xcolor}
\usepackage{mathrsfs}
\usepackage{dsfont}

\newtheorem{theorem}{Theorem}[section]
\newtheorem{lemma}[theorem]{Lemma}
\newtheorem{proposition}[theorem]{Proposition} 
\newtheorem{corollary}[theorem]{Corollary}
\newtheorem{remark}[theorem]{Remark}
\newtheorem{definition}[theorem]{Definition}

\title{\LARGE \bf
On exponential stabilization of $N$-level quantum angular momentum systems
}
\date{\vspace{-5ex}}

\author{Weichao Liang, Nina H. Amini, and Paolo Mason
\thanks{All authors are with Laboratoire des Signaux et Syst\`emes, CNRS - CentraleSup\'elec - Univ. Paris-Sud, Universit\'e Paris-Saclay, 3, rue Joliot Curie, 91192, Gif-sur-Yvette, France. {\tt\small [first name].[family name]@l2s.centralesupelec.fr}}%
}

\begin{document}

\maketitle

\begin{abstract}
In this paper, we consider the feedback stabilization problem for $N$-level quantum angular momentum systems undergoing continuous-time measurements. By using stochastic and geometric control tools, we provide sufficient conditions on the feedback control law ensuring almost sure exponential convergence to a predetermined eigenstate of the measurement operator. In order to achieve these results, we establish general features of quantum trajectories which are of interest by themselves. We illustrate the results by designing a class of feedback control laws satisfying the above-mentioned conditions and finally we demonstrate the effectiveness of our methodology through numerical simulations for three-level quantum angular momentum systems.
\end{abstract}

\section{Introduction}
The evolution of an open quantum system undergoing indirect continuous-time measurements is described by the so-called quantum stochastic master equation, which has  been derived by Belavkin in quantum filtering theory \cite{belavkin1989nondemolition}. The quantum filtering theory, relying on quantum stochastic calculus and quantum probability theory (developed by Hudson and Parthasarathy~\cite{hudson1984quantum}) plays an important role in quantum optics and computation. The initial concepts of quantum filtering have been developed in the 1960s by Davies~\cite{davies1969quantum,davies1976quantum} and extended by Belavkin in the 1980s~\cite{belavkin1983theory,belavkin1989nondemolition,belavkin1995quantum,belavkin1992quantum}. For a modern treatment of quantum filtering, we refer to~\cite{bouten2007introduction,van2005feedback}. 

A quantum stochastic master equation (or quantum filtering equation) is composed of a deterministic part and a stochastic part. The deterministic part, which corresponds to the average dynamics, is given by the well known Lindblad operator. The stochastic part represents the back-action effect of continuous-time measurements. The solutions of this equation are called quantum trajectories and their properties have been studied in~\cite{mirrahimi2007stabilizing, pellegrini2008existence}.

Quantum measurement-based feedback control, as a branch of stochastic control has been first developed by Belavkin in~\cite{belavkin1983theory}.  This field has attracted the interest of many theoretical and experimental researchers mainly starting from the early 2000s, yielding fundamental results~\cite{van2005feedback, armen2002adaptive, mirrahimi2007stabilizing, tsumura2008global,ahn2002continuous,yamamoto2007feedback,mabuchi2005principles}. In particular, theoretical studies carried out in~\cite{mirrahimi2007stabilizing,dotsenko2009quantum,mirrahimi2009feedback,amini2011design,amini2013feedback} lead to the first experimental implementation of real-time quantum measurement-based feedback control in~\cite{sayrin2011real}. 

In~\cite{bouten2008separation}, the authors established a quantum separation principle. Similarly to the classical separation principle,  this result allows to interpret the control problem as a state-based feedback control problem for the filter (the best estimate, i.e., the conditional state), without caring of the actual quantum state. This motivates the state-based feedback design for the quantum filtering equation based on the knowledge of the initial state.
In this context, stabilization of quantum filters towards pure states (i.e., the preparation of pure states) has major impact in developing new quantum technologies.  According to~\cite{abe2008analysis}, the stochastic part of the quantum filtering equation, unlike the deterministic one, contributes to increase the purity of the quantum state. Moreover, if we turn off the control acting on the quantum system, the measurement induces a collapse of the quantum state towards either one of the eigenstates of the measurement operator, a phenomenon known as quantum state reduction~\cite{adler2001martingale, van2005feedback, mirrahimi2007stabilizing, sarlette2017deterministic}. Thus, combining the continuous measurement with the feedback control may provide an effective strategy for preparing a selected target eigenstate  in practice.

In~\cite{van2005feedback}, the authors design for the first time a quantum feedback controller that globally stabilizes a quantum spin-$\frac{1}{2}$ system (which is a special case of quantum angular momentum systems) towards an eigenstate of $\sigma_z$ in the presence of imperfect measurement. This feedback controller has been designed by looking numerically for an appropriate global Lyapunov function. Then, in ~\cite{mirrahimi2007stabilizing}, by analyzing the stochastic flow and by using stochastic Lyapunov techniques, the authors constructed a switching feedback controller which globally stabilizes the $N$-level quantum angular momentum system, in the presence of imperfect measurement, to the target eigenstate. A continuous version of this feedback controller has been proposed in~\cite{tsumura2008global}. The essential ideas in~\cite{van2005feedback, tsumura2008global} for constructing the continuous feedback controller remain the same: the controllers consist of two parts, the first one contributing to the local convergence to the target eigenstate, and the second one driving the system away from the antipodal eigenstates. Also, in~\cite{cardona2018exponential}, the authors have proven by simple Lyapunov arguments the stochastic exponential stabilizability for spin-$\frac{1}{2}$ systems  by applying a proportional output feedback.

The main contribution of this paper is the derivation of  some general conditions on the feedback law enforcing the exponential convergence towards the target state. These conditions are obtained mainly by studying the asymptotic behavior of quantum trajectories. Roughly speaking, under such conditions, and making use of the support theorem and other  classical stochastic tools, we show that any neighborhood of the target state may be approached with non-zero probability starting from any initial state. The exponential convergence towards the target state is then obtained via Lyapunov arguments. As demonstration of the general result, explicit parametrized stabilizing feedback laws are exhibited. In addition to the main result, we show the exponential convergence of the system with zero control towards the set of eigenstates of the measurement operator (quantum state reduction with exponential rate).  Note that to obtain our main results, some preliminary results on  the asymptotic behavior of quantum trajectories associated with the considered system were needed.  We believe that these results are  significants by themselves. We point out that preliminary results for two-level angular momentum systems were provided in~\cite{liang2018exponential}.
\paragraph{Notations}
The imaginary unit is denoted by $i$. We take $\mathds{1}$ as the indicator function. We denote the conjugate transpose of a matrix $A$ by $A^*.$ The function $\mathrm{Tr}(A)$ corresponds to the trace of a square matrix $A.$ The commutator of two square matrices $A$ and $B$ is denoted by $[A,B]:=AB-BA.$

\section{System description}
Consider a filtered probability space $(\Omega,\mathcal{F},(\mathcal{F}_t),\mathbb{P})$. Let $W_t$ be the 1-dimensional standard Wiener process and assume that $\mathcal{F}_t$ is the natural filtration of the process $W_t$. The dynamics of a $N$-level quantum angular momentum system is given by the following matrix-valued stochastic differential equation~\cite{belavkin1989nondemolition,bouten2007introduction,van2005feedback}:
\begin{equation}
d\rho_t=F(\rho_t)dt+\sqrt{\eta}G(\rho_t)dW_t,
\label{ND SME}
\end{equation}
where
\begin{itemize}
\item The quantum state is described by the density operator $\rho$, which belongs to the compact space $\mathcal{S}:=\{\rho\in\mathbb{C}^{N \times N}|\,\rho=\rho^*,\mathrm{Tr}(\rho)=1,\rho\geq0\}$, 
\item the drift term is given by $$F(\rho):=-i\omega[J_z,\rho]+M\left(J_z\rho J_z-\frac12J^2_z\rho-\frac12\rho J^2_z\right)-iu_t[J_y,\rho],$$
and the diffusion term is given by $G(\rho):=\sqrt{M}(J_z\rho+\rho J_z-2\mathrm{Tr}(J_z\rho)\rho),$
\item $u_t:=u(\rho_t)$ denotes the feedback law, 
\item 
$J_z$ is the (self-adjoint) angular momentum along the axis $z$, and it is defined by 
\begin{equation*}
J_z e_n=(J-n)e_n,\quad n\in\{0,\dots,2J\}, 
\end{equation*}
where $J=\frac{N-1}{2}$ represents the fixed angular momentum and $\{e_0,\dots,e_{2J}\}$ corresponds to an orthonormal basis of $\mathbb C^N.$ With respect to this basis, the matrix form of $J_z$ is given by
\begin{equation*}
J_z=
\begin{bmatrix}
J &&&&  \\
& J-1&&&\\
&&\ddots&&\\
&&&-J+1&\\
&&&&-J
\end{bmatrix},
\end{equation*}
\item 
$J_y$ is the (self-adjoint) angular momentum along the axis $y$, and it is defined by 
\begin{equation*}
\begin{split}
J_ye_n=-ic_{n}e_{n-1}+ic_{n+1}e_{n+1},\quad n\in\{0,\dots,2J\},
\end{split}
\end{equation*}
where $c_m=\frac12\sqrt{(2J+1-m)m}$. The matrix form of $J_y$ is given by
\begin{equation*}
J_y=
\begin{bmatrix}
0&-ic_1 &&&\\
ic_1&0&-ic_2&&\\
&\ddots&\ddots&\ddots&\\
&&ic_{2J-1}&0&-ic_{2J}\\
&&&ic_{2J}&0
\end{bmatrix},
\end{equation*}
\item $\eta\in(0,1]$ measures the efficiency of the photon-detectors,  $M>0$ is the strength of the interaction between the light and the atoms, and $\omega \geq 0$ is a parameter characterizing the free Hamiltonian. 
\end{itemize} 
If the feedback $u$ is in $\mathcal{C}^1(\mathcal{S},\mathbb{R})$, the existence and uniqueness of the solution of~\eqref{ND SME} as well as the strong Markov property of the solution are ensured by the results established in~\cite{mirrahimi2007stabilizing}. 

\section{Basic stochastic tools}
In this section, we will introduce some basic definitions and classical results which are fundamental for the rest of the paper.
\paragraph{Infinitesimal generator and It\^o's formula}
Given a stochastic differential equation $dq_t=f(q_t)dt+g(q_t)dW_t$, where $q_t$ takes values in $Q\subset \mathbb{R}^p,$ the infinitesimal generator  is the operator $\mathscr{L}$ acting on twice continuously differentiable functions $V: Q \times \mathbb{R}_+ \rightarrow \mathbb{R}$ in the following way
\begin{equation*}
\mathscr{L}V(q,t):=\frac{\partial V(q,t)}{\partial t}+\sum_{i=1}^p\frac{\partial V(q,t)}{\partial q_i}f_i(q)+\frac12 \sum_{i,j=1}^p\frac{\partial^2 V(q,t)}{\partial q_i\partial q_j}g_i(q)g_j(q).
\end{equation*}
It\^o's formula describes the variation of the function $V$ along solutions of the stochastic differential equation and is given as follows
\begin{equation*}
dV(q,t) = \mathscr{L}V(q,t)dt+\sum_{i=1}^p\frac{\partial V(q,t)}{\partial q_i}g_i(q)dW_t.
\end{equation*}
From now on, the operator $\mathscr{L}$ is associated with the equation~\eqref{ND SME}.
 
\paragraph{Stochastic stability}
We introduce some notions of stochastic stability needed throughout the paper by adapting classical notions (see e.g.~\cite{mao2007stochastic, khasminskii2011stochastic}) to our setting. In order to provide them, we first present the definition of Bures distance~\cite{bengtsson2017geometry}.
\begin{definition}
The Bures distance between two quantum states $\rho_a$ and $\rho_b$ in $\mathcal{S}$ is defined as
\begin{equation*}
d_B(\rho_a,\rho_b) := \sqrt{2-2\mathrm{Tr}\left( \sqrt{\sqrt{\rho_b}\rho_a\sqrt{\rho_b}} \right)}.
\end{equation*}
In particular, the Bures distance between a quantum state $\rho_a \in \mathcal{S}$ and a pure state $\boldsymbol\rho_{n}:=e_ne^*_n$ with $n\in\{0,\dots,2J\},$ is given by
\begin{equation*}
d_B(\rho_a,\boldsymbol\rho_{n})= \sqrt{2-2\sqrt{\mathrm{Tr}(\rho_a\boldsymbol\rho_{n})}}.
\end{equation*}
Also, the Bures distance between a quantum state $\rho_a$ and a set $E \subseteq \mathcal{S}$ is defined as 
\begin{equation*}
d_B(\rho_a, E) = \min_{\rho \in E} d_B(\rho_a,\rho).
\end{equation*}
\end{definition}
Given $E \subseteq \mathcal{S}$ and $r>0$, we define the neighborhood $B_r(E)$ of $E$ as
\begin{equation*}
B_r(E) = \{\rho \in \mathcal{S}|\, d_B(\rho,E) < r\}.
\end{equation*}
\begin{definition}
Let $\bar E$ be an invariant set of system~\eqref{ND SME}, then $\bar E$ is said to be
\begin{enumerate}
\item[1.] 
\emph{locally stable in probability}, if for every $\varepsilon \in (0,1)$ and for every $r >0$, there exists  $\delta = \delta(\varepsilon,r)$ such that,
\begin{equation*}
\mathbb{P} \left( \rho_t \in B_r (\bar E) \text{ for } t \geq 0 \right) \geq 1-\varepsilon,
\end{equation*}
whenever $\rho_0 \in B_{\delta} (\bar E)$.

\item[2.]
\emph{exponentially stable in mean}, if for some positive constants $\alpha$ and $\beta$,
\begin{equation*}
\mathbb{E}(d_B(\rho_t,\bar E)) \leq \alpha \,d_B(\rho_0,\bar E)e^{-\beta t},
\end{equation*} 
whenever $\rho_0 \in \mathcal{S}$. The smallest value $-\beta$ for which the above inequality is satisfied is called the \emph{average Lyapunov exponent}. 

\item[3.]
\emph{almost surely exponentially stable}, if
\begin{equation*}
\limsup_{t \rightarrow \infty} \frac{1}{t} \log d_B(\rho_t,\bar E) < 0, \quad a.s.
\end{equation*}
whenever $\rho_0 \in \mathcal{S}$. The left-hand side of the above inequality is called the \emph{sample Lyapunov exponent} of the solution. 
\end{enumerate}
\end{definition}
Note that any equilibrium $\bar\rho$ of~\eqref{ND SME}, that is any quantum state satisfying $F(\bar\rho)=G(\bar\rho)=0$, is a special case of invariant set.
\paragraph{Stratonovich equation and Support theorem}
Any stochastic differential equation in It\^o form in $\mathbb R^K$
\begin{equation*}
dx_t=\widehat X_0(x_t)dt+\sum^n_{k=1}\widehat X_k(x_t)dW^k_t, \quad x_0 = x,
\end{equation*}
can be written in the following Stratonovich form~\cite{rogers2000diffusions2}
\begin{equation*}
dx_t = X_0(x_t)dt+\sum^n_{k=1}X_k(x_t) \circ dW^k_t, \quad x_0 = x,
\end{equation*}
where 
$X_0(x)=\widehat X_0(x)-\frac{1}{2}\sum^K_{l=1}\sum^n_{k=1}\frac{\partial \widehat X_k}{\partial x_l}(x)(\widehat X_k)_l(x)$, $(\widehat X_k)_l$ denoting the component $l$ of the vector $\widehat X_k,$ and $X_k(x)=\widehat X_k(x)$ for $k\neq 0$.

\medskip
The following classical theorem relates the solutions of a stochastic differential equation with those of an associated deterministic one.
\begin{theorem}[Support theorem~\cite{stroock1972support}]
Let $X_0(t,x)$ be a bounded measurable function, uniformly Lipschitz continuous in $x$ and $X_k(t,x)$  be continuously differentiable in $t$ and twice continuously differentiable in $x$, with bounded derivatives, for $k\neq 0.$ Consider the Stratonovich equation
\begin{equation*}
dx_t = X_0(t,x_t)dt+\sum^n_{k=1}X_k(t,x_t) \circ dW^k_t, \quad x_0 = x.
\end{equation*}
Let $\mathbb{P}_x$ be the probability law of the solution $x_t$ starting at $x$. Consider in addition the associated deterministic control system
\begin{equation*}
\frac{d}{dt}x_{v}(t) = X_0(t,x_{v}(t))+\sum^n_{k=1}X_k(t,x_{v}(t))v^k(t), \quad x_v(0) = x.
\end{equation*}
with $v^k \in \mathcal{V}$, where $\mathcal{V}$ is the set of all piecewise constant functions from $\mathbb{R}_+$ to $\mathbb{R}$. Now we define $\mathcal{W}_x$ as the set of all continuous paths from $\mathbb{R}_+$ to $\mathbb R^K$ starting at $x$, equipped with the topology of uniform convergence on compact sets, and $\mathcal{I}_x$ as the smallest closed subset of $\mathcal{W}_x$ such that $\mathbb{P}_x(x_{\cdot} \in \mathcal{I}_x)=1$. Then,
\begin{equation*}
\mathcal{I}_x = \overline{ \{ x_{v}(\cdot)\in\mathcal{W}_x|\, v \in \mathcal{V}^n\} } \subset \mathcal{W}_x.
\end{equation*}
\label{Support thm}
\end{theorem}

\section{Preliminary results}
Our aim here is to establish some basic properties of the quantum trajectories corresponding to Equation~\eqref{ND SME}. This section is instrumental  in order to prove our main results. 

Denote the projection of $\rho$ onto the eigenstate $\boldsymbol\rho_{k}$ as $\rho_{k,k}:=\mathrm{Tr}(\rho\boldsymbol\rho_{k}).$ In the following we state two lemmas inspired by analogous results established in~\cite{khasminskii2011stochastic, mao2007stochastic}. 
\begin{lemma}
Assume $u\equiv 0$.  If $\rho_{k,k}(0)=0$ for some $k\in \{0,\dots,2J\},$ then $\mathbb{P}( \rho_{k,k}(t)=0, \forall\, t\geq 0 )=1,$ i.e., the set $\{\rho\in\mathcal S|\, \rho_{k,k}=0\}$ is a.s. invariant for Equation~\eqref{ND SME}. Otherwise, if the initial state satisfies $\rho_{k,k}(0)\neq 0$, then
$
\mathbb{P}( \rho_{k,k}(t)\neq0, \forall\, t\geq 0 )=1.$ 
\label{Never reach without feedback} 
\end{lemma}
\proof
For $u\equiv 0,$ the dynamics of $\rho_{k,k}$ is given by 
\begin{equation*}
d\rho_{k,k}(t)=\sqrt{\eta}(G(\rho_t))_{k,k}dW_t=2\sqrt{\eta M}(J-k-\mathrm{Tr}(J_z\rho_t))\rho_{k,k}(t)dW_t.
\end{equation*}
In particular $|\sqrt{\eta}(G(\rho_t))_{k,k}|\leq R\rho_{k,k}(t),$ for some $R>0,$ yielding the first part of the lemma.

Let us now prove the second part of the lemma. Assume that $\rho_{k,k}(0)>0$ and $\mathbb{P}( \rho_{k,k}(t)\neq0, \forall\, t\geq 0 )<1.$ In particular $\mathbb{P}(\tau<\infty)>0$, where $\tau:=\inf\{t\geq 0|\,\rho_{k,k}(t)=0\}$.
Let $T$ be sufficiently large so that $\mathbb{P}(\tau\leq T)>0$. Now, let $\varepsilon\in (0,\rho_{k,k}(0))$, and consider any $\mathcal{C}^2$ function $V$ defined on $\mathcal{S}$ such that
\begin{equation*}
V(\rho)=\frac{1}{\rho_{k,k}},\qquad \mbox{if }\rho_{k,k}>\varepsilon.
\end{equation*}
Then we have $\mathscr{L}V(\rho)=\rho^{-3}_{k,k}(\sqrt{\eta}G(\rho))^2_{k,k}\leq R^2V(\rho)$ if $\rho_{k,k}>\varepsilon$. We further define the time-dependent function 
$
f(\rho,t)=e^{-R^2t}V(\rho),
$
whose infinitesimal generator is given by
$
\mathscr{L}f(\rho,t)=e^{-R^2t}\left(-R^2V(\rho)+\mathscr{L}V(\rho)\right)\leq 0
$
if $\rho_{k,k}>\varepsilon.$
Now, define the stopping time $\tau_{\varepsilon}:=\inf\{t \geq0|\,\rho_{k,k}(t) \notin (\varepsilon,1)\}$. By It\^o's formula, we have
\begin{equation*}
\mathbb{E}(f(\rho_{\tau_{\varepsilon}\wedge T},\tau_{\varepsilon}\wedge T))=V_0+\mathbb{E}\left(\int^{\tau_{\varepsilon}\wedge T}_{0}\mathscr{L}f(\rho_s,s)ds\right)\leq V_0=\frac1{\rho_{k,k}(0)}.
\end{equation*}
Since $\tau\geq \tau_{\varepsilon}$ we deduce that, conditioning to the event $\{\tau\leq T\}$, $f(\rho_{\tau_{\varepsilon}\wedge T},\tau_{\varepsilon}\wedge T) = f(\rho_{\tau_{\varepsilon}},\tau_{\varepsilon}) = e^{-R^2T}\varepsilon^{-1}$,
which implies
\begin{equation*}
\mathbb{E}\left(e^{-R^2T}\varepsilon^{-1}\mathds{1}_{\{\tau\leq T\}}\right) = \mathbb{E}\left(f(\rho_{\tau_{\varepsilon}},\tau_{\varepsilon})\mathds{1}_{\{\tau\leq T\}}\right)\leq \mathbb{E}(f(\rho_{\tau_{\varepsilon}\wedge T},\tau_{\varepsilon}\wedge T))\leq \frac1{\rho_{k,k}(0)}.
\end{equation*}
Thus,
$
\mathbb{P}(\tau\leq T)=\mathbb{E}\left(\mathds{1}_{\{\tau\leq T\}}\right)\leq \varepsilon e^{R^2T}/\rho_{k,k}(0).
$
Letting $\varepsilon$ tend to $0$, we get $\mathbb{P}(\tau\leq T)=0$ which gives a contradiction. The proof is then complete. \hfill$\square$
\begin{lemma}
Let $n \in \{0,\dots,2J\}.$ Assume that the initial state satisfies $\rho_0 \neq \boldsymbol \rho_n,$ $u\in\mathcal{C}^1(\mathcal{S}\setminus\boldsymbol \rho_n,\mathbb{R})$ and $u(\rho)\leq C\sqrt{1-\rho_{n,n}}$ for some $C>0$ . Then
$
\mathbb{P}( \rho_t \neq \boldsymbol \rho_n, \forall\, t\geq 0 )=1.
$
\label{Never reach with feedback}
\end{lemma}
\proof
Given $\varepsilon>0,$ we consider any $\mathcal{C}^2$ function on $\mathcal{S}$ such that
\begin{equation*}
V(\rho)=\frac{1}{1-\rho_{n,n}},\qquad \mbox{if } \rho_{n,n}<1-\varepsilon.
\end{equation*}
We find
\begin{equation*}
\mathscr{L}V(\rho)=-\frac{u(\rho)\mathrm{Tr}(i[J_y,\rho]\boldsymbol \rho_n)}{(1-\rho_{n,n})^2}+\frac{4\eta M [(J-n-\mathrm{Tr}(J_z\rho))\rho_{n,n}]^2}{(1-\rho_{n,n})^3},
\end{equation*}
whenever $ \rho_{n,n}<1-\varepsilon$. Since 
\begin{equation*}
\begin{split}
\mathrm{Tr}(i[J_y,\rho]\boldsymbol \rho_n)&=2c_{n+1}\mathbf{Re}\{\rho_{n,n+1}\}-2c_n\mathbf{Re}\{\rho_{n,n-1}\}\\
&\leq 2(c_{n+1}+c_n)\sqrt{\rho_{n,n}(1-\rho_{n,n})}
\end{split}
\end{equation*} 
and $u(\rho)\leq C\sqrt{1-\rho_{n,n}}$, we have
$
|u(\rho)\mathrm{Tr}(i[J_y,\rho]\boldsymbol \rho_n)|\leq 2C(c_{n+1}+c_n)(1-\rho_{n,n}).
$
Also, as we have $|J-n-\mathrm{Tr}(J_z\rho)|\leq 2J(1-\rho_{n,n}),$ we get
$\mathscr{L}V(\rho)\leq KV(\rho),$ with $K = 2C(c_{n+1}+c_n)+16J^2\eta M.$
To conclude the proof, one just applies the same arguments as in the previous lemma. \hfill$\square$

\medskip
Consider the observation process of the system $y_t,$ whose dynamics satisfies $dy_t=dW_t+2\sqrt{\eta M}\mathrm{Tr}(J_z\rho_t)dt$.
By Girsanov's theorem~\cite{protter2004stochastic}, the process $y_t$ is a standard Wiener process under a new probability measure $\mathbb{Q}$ equivalent to $\mathbb P$. Denote by $\mathcal{F}^y_t:=\sigma(y_s,0\leq s\leq t)$ the $\sigma$-field generated by the observation process up to time $t$. Then by applying the classical stochastic filtering theory~\cite{xiong2008introduction}, the Zakai equation associated with Equation~\eqref{ND SME} takes the following linear form
\begin{equation}
d\tilde{\rho}_t=F(\tilde{\rho}_t)dt+\sqrt{\eta}\widetilde{G}(\tilde{\rho}_t)dy_t,
\label{Zakai ND}
\end{equation} 
where $\tilde{\rho}_t=\tilde{\rho}^*_t\geq0$, 
$F(\tilde{\rho})$ is defined as in~\eqref{ND SME}, and $\widetilde{G}(\tilde{\rho}):=\sqrt{M}(J_z\tilde{\rho}_t+\tilde{\rho}_tJ_z).$
The equation~\eqref{Zakai ND} has a unique strong solution~\cite{xiong2008introduction,protter2004stochastic}, and the solutions of the equations~\eqref{ND SME} and~\eqref{Zakai ND} satisfy the relation
\begin{equation}
\rho_t=\tilde{\rho_t}/\mathrm{Tr}(\tilde{\rho_t}),
\label{Relation Zakai FFk}
\end{equation}
which can be verified easily by applying It\^o's formula. In the following lemma, we adapt~\cite[Lemma 3.2]{mirrahimi2007stabilizing} to the case of positive-definite matrices.
\begin{lemma}
The set of positive-definite matrices is a.s. invariant for~\eqref{ND SME}. More in general, the rank of $\rho_t$ is a.s. non-decreasing.
\label{Pos def invariant}
\end{lemma}
\proof
The initial state of~\eqref{Zakai ND} with respect to the basis of its eigenstates is given by $\tilde{\rho}_0=\sum_i\tilde{\lambda}_i\tilde{\psi}_i\tilde{\psi}^*_i$, where $\tilde{\rho}_0\tilde{\psi}_i=\tilde{\lambda}_i\tilde{\psi}_i$ for $i\in\{0,\dots,2J\}$. If $\rho_0>0$, due to the relation~\eqref{Relation Zakai FFk}, we have $\tilde{\rho}_0>0$, thus $\tilde{\lambda}_i>0$ for all $i$. Extend the probability space by defining $\mathcal{F}^{y,\widetilde{W}}_t:=\sigma(y_s,\widetilde{W}_s,0\leq s\leq t)$, where $\widetilde{W}_t$ is a Brownian motion independent of $y_t$. Set $B_t:=\sqrt{\eta}y_t+\sqrt{1-\eta}\widetilde{W}_t$, whose quadratic variation satisfies $\left\langle B_t,B_t\right\rangle=t$. Following~\cite[Lemma 3.2]{mirrahimi2007stabilizing}, we consider the equations
\begin{equation*}
\begin{split}
&d\rho^i_t = F(\rho^i_t)dt+\widetilde{G}(\rho^i_t)\sqrt{\eta}dy_t+\widetilde{G}(\rho^i_t)\sqrt{1-\eta}\,d\widetilde{W}_t,\quad \rho^i_0=\tilde{\psi}_i\tilde{\psi}^*_i,\\
&d\tilde{\psi}_i(t) = (i\omega J_z-iu_tJ_y-M/2J^2_z)\tilde{\psi}_i(t)dt+\sqrt{M}J_z\tilde{\psi}_i(t)dB_t,\quad \tilde{\psi}_i(0)=\tilde{\psi}_i, 
\end{split}
\end{equation*}
where $\tilde{\psi}_i(t)\in\mathbb{C}^N.$ The solutions of the equations above satisfy $\rho^i_t = \tilde{\psi}_i(t)\tilde{\psi}^*_i(t)$ by It\^ o's formula. In virtue of~\cite[Theorem~5.48]{protter2004stochastic}, for all $t\geq0$, there exists an almost surely invertible random matrix $U_t$ such that $\tilde{\psi}_i(t)=U_t\tilde{\psi}_i$.

Let $\rho'_t=\sum_i\tilde{\lambda}_i\rho^i_t$, so that in particular $\rho'_0 = \tilde{\rho}_0$ and $\rho'_t=U_t\tilde{\rho}_0U^*_t$. Due to the linearity of $F(\cdot)$ and $\widetilde{G}(\cdot)$, the stochastic Fubini theorem~\cite[Lemma 5.4]{xiong2008introduction} and the It\^o's isometry, 
\begin{equation*}
\mathbb{E}(\rho'_t|\mathcal{F}^y_t) = \rho'_0 + \int^t_0 F(\mathbb{E}(\rho'_s|\mathcal{F}^y_t))ds+\int^t_0\widetilde{G}(\mathbb{E}(\rho'_s|\mathcal{F}^y_t))\sqrt{\eta}dy_s.
\end{equation*}
By the uniqueness in law~\cite[Proposition 9.1.4]{revuz2013continuous} of the solution of the equation~\eqref{Zakai ND}, the laws of $\tilde{\rho}_t$ and $\mathbb{E}(\rho'_t|\mathcal{F}^y_t)= \mathbb{E}(U_t\tilde{\rho}_0U^*_t|\mathcal{F}^y_t)$ are equal for all $t\geq 0$. 

By what precedes $\rho_0>0$ implies $\rho'_t>0$ a.s. which in turn yields $\rho_t = \tilde{\rho}_t/\mathrm{Tr}(\tilde{\rho_t})>0$ a.s. We have thus proved that the set of positive-definite matrices is a.s. invariant for~\eqref{ND SME}.

Let us now consider the general case in which $\rho_0$ is not necessarily full rank. We have 
\begin{equation}
\mathrm{rank}(\rho'_t)=\mathrm{rank}(U_t\tilde{\rho}_0U^*_t)=\mathrm{rank}(\tilde{\rho}_0)=\mathrm{rank}(\rho_0),\quad a.s.
\label{eq:rank}
\end{equation}
Note that the kernel of any positive semi-definite matrix $\hat\rho\in\mathbb C^{N\times N}$ coincides with the space $\{\psi\in \mathbb C^N|\psi^*\hat\rho\psi=0\}$, and that for almost every path $\rho'_t(\omega)$
\begin{equation*}
\{\psi\in \mathbb C^N|\mathbb{E}(\psi^*\rho'_t\psi|\mathcal{F}^y_t)=0\}\subseteq \{\psi\in \mathbb C^N|\psi^*\rho'_t(\omega)\psi=0\}.
\end{equation*}
This implies 
$
\mathrm{rank}(\tilde\rho_t)\geq \mathrm{rank}(\rho'_t)=\mathrm{rank}(\rho_0)
$
for any $t\geq 0$ almost surely, which concludes the proof.\hfill$\square$

\begin{lemma}
If $\eta=1$, then the boundary of the state space 
\begin{equation*}
\partial\mathcal{S}:=\{\rho\in\mathbb{C}^{N\times N}|\,\rho=\rho^*,\mathrm{Tr}(\rho)=1,\det(\rho)=0\}
\end{equation*}
is a.s. invariant for~\eqref{ND SME}.
\label{boundary}
\end{lemma}
\proof
Based on the proof of Lemma~\ref{Pos def invariant}, if $\eta=1$, we have $B_t=y_t$ which implies $\tilde{\rho}_t=\rho'_t$. Then by applying the relation~\eqref{eq:rank}, we get the conclusion. \hfill$\square$

\medskip

The Stratonovich form of Equation~\eqref{ND SME} is given by
\begin{equation}
d\rho_t=\widehat{F}(\rho_t)dt+\sqrt{\eta}G(\rho_{t})\circ dW_t,
\label{ND strat}
\end{equation}
where 
\begin{equation*}
\begin{split}
\widehat{F}(\rho):=&-i\omega[J_z,\rho]+M\left((1-\eta)J_z\rho J_z-\frac{1+\eta}{2}(J^2_z\rho+\rho J^2_z)+2\eta \mathrm{Tr}(J^2_z\rho)\rho\right)\\
&+2\eta M\mathrm{Tr}(J_z \rho)(J_z\rho+\rho J_z-2\mathrm{Tr}[J_z\rho]\rho)-iu(\rho)[J_y,\rho],\\
\end{split}
\end{equation*}
and $G$ is defined as in~\eqref{ND SME}. The corresponding deterministic control system is given by
\begin{equation}
\dot{\rho}_{v}(t)=\widehat{F}(\rho_{v}(t))+\sqrt{\eta}G(\rho_{v}(t))v(t),\quad \rho_{v}(0)=\rho_0,
\label{ODE}
\end{equation}
where $v(t)\in\mathcal{V}$. By the support theorem (Theorem~\ref{Support thm}), the set $\mathcal S$ is positively invariant for Equation~\eqref{ODE}.

In the following, we state some preliminary results that will be applied to our stabilization problem in the following sections. For this purpose, we fix a target state $\boldsymbol \rho_{\bar n}$ for some  $\bar n\in\{0,\dots,2J\}.$
\begin{proposition}
Suppose $\eta\in(0,1)$ and $u\in\mathcal{C}^1(\mathcal{S}\setminus\boldsymbol \rho_{\bar n},\mathbb{R})$. Assume that $\nabla u\cdot G(\rho_0)\neq 0$ or $\nabla u\cdot \widehat F(\rho_0)\neq 0$ for any $\rho_0\in\{\rho\in\mathcal S\setminus{\boldsymbol \rho_{\bar n}}| \,\, \rho_{k,k}=0 \,\ \textrm{for some $k,$ and} \ u(\rho)=0\}.$ Then for any initial condition $\rho_0\in\{\rho\in\mathcal S\setminus{\boldsymbol \rho_{\bar n}}| \,\, \rho_{k,k}=0 \,\ \textrm{for some $k$}\}$ and $\varepsilon>0,$ there exists at most one trajectory $\rho_{v}(t)$ of \eqref{ODE} starting from $\rho_0$ which lies in $\partial\mathcal{S}$ for $t$ in $[0,\varepsilon]$. For any other initial state $\rho_0\in\partial\mathcal{S}\setminus{\boldsymbol \rho_{\bar n}}$ and $v\in\mathcal{V},$ $\rho_{v}(t)>0$ for $t>0$.
\label{Exits boundary ODE}
\end{proposition}
\proof
Define $Z_1(t):=\textrm{Span}\{e_k|\, (\rho_{v}(t))_{k,k}=0\}$ and $Z_2(t)$ the eigenspace corresponding to the eigenvalue $0$ of $\rho_{v}(t).$ By definition, $Z_1(t)\subseteq Z_2(t)$ for all $t\geq 0.$ Since all the subspaces which are invariant by $J_z$ take the form $\textrm{Span}\{e_{k_1},\dots,e_{k_h} \}$ for $\{k_1,\dots,k_h\}\subseteq\{0,\dots,2J\},$ we deduce that $Z_1(t)$  is the largest subspace of $Z_2(t)$ invariant by $J_z.$ 

Denote by $\lambda_{k}(t)$ and $\psi_{k}(t)$ for $k\in\{0,\dots,2J\}$ the eigenvalues and eigenvectors of $\rho_{v}(t)$, where, without loss of generality, we assume $\lambda_k(t)\in\mathcal{C}^1$ since  $\rho_{v}(t)\in\mathcal{C}^1$ (\cite[Theorem 2.6.8]{kato1976perturbation}). In addition, we suppose that the eigenvectors $\psi_{k}(t)$ form an orthonormal basis of $\mathbb C^N.$

Let $\psi_k(t)\in Z_2(t)$ for $t\in[0,\varepsilon].$   In order to provide an expression of  the derivative for the eigenvalue $\lambda_{k}$ along the path, we observe that
\begin{equation}
\begin{split}\label{derivativelambda}
\frac1t (\lambda_k(t+\delta)-\lambda_k(t))   
= \frac1{\psi_k^*(t+\delta)\psi_k(t)} \left(\psi_k^*(t+\delta) \frac{\rho_{v}(t+\delta)-\rho_{v}(t)}{t}\psi_k(t)\right).
\end{split}
\end{equation}
Since $\psi_k$ is a unit vector, then by compactness, we can extract a sequence $\delta_n\searrow 0$ such that  $\psi_k(t+\delta_n)$ converges to an eigenvector $\psi_k(t)$ of $\rho_{v}(t)$. By passing to the limit on the left-hand and right-hand sides of Equation~\eqref{derivativelambda}, we get $\dot{\lambda}_k(t) = \psi^*_{k}(t)\dot{\rho}_{v}(t)\psi_{k}(t)=M(1-\eta)\psi^*_k(t)J_z\rho_{v}(t)J_z\psi_{k}(t)$. 

If $\psi_k(t)\notin Z_1(t)$ then $J_z\psi_k(t)\notin Z_2(t),$ since otherwise $Z_1(t)$ would not be the largest subspace  invariant by $J_z$ contained in $Z_2(t).$ Thus $\dot{\lambda}_k(t)>0,$ which implies $\lambda_k(s)>0$ for any $s-t>0$ sufficiently small. We deduce that $\textrm{dim}\,Z_2(s)\leq \textrm{dim}\,Z_1(t)$.  Moreover, by continuity of $\rho_{v}(t),$ we have $Z_1(s) \subseteq Z_1(t),$ for any $s-t>0$ sufficiently small. Now we consider the case where $Z_1(t)\neq 0$ for $t\geq 0.$ In this case, we have two possibilities:  either $u(\rho_{v}(\cdot))\equiv 0$ on $[0,\varepsilon]$ for some $\varepsilon>0;$ or 
$u(\rho_{v}(t))\neq 0$ for arbitrarily small $t>0.$ Note that under the assumptions of the proposition there exists at most one $v$ such that $u(\rho_{v}(\cdot))\equiv 0.$  It is therefore enough to show that, for the second possibility, $\rho_{v}(t)$ belongs to the interior of $\mathcal S$ for all $t>0$. For this purpose, we first show that for all $t>0$ such that $u(\rho_{v}(t))\neq 0$ and $Z_1(t)\neq 0$, there exists $s-t>0$ arbitrarily small such that $u(\rho_{v}(s))\neq 0$ and $Z_1(s)\subsetneqq Z_1(t)$. 

Let us pick $k$ such that $e_{k}\in Z_1(t),$ and at least one between $e_{k-1}$ and $e_{k+1}$ is not contained in $Z_1(t)$\footnote{If $k=0,$ the condition is replaced by $e_1\notin Z_1(t)$ while if $k=2J,$ we assume $e_{2J-1}\notin Z_1(t)$.}. We now show by contradiction that $e_{k}\notin Z_1(s)$ for some $s-t>0$ arbitrarily small. We assume that $e_k\in Z_1(\tau)$ for $\tau\in [t, t+\varepsilon]$, with $\varepsilon>0.$ By setting $q^n({\tau}):=\rho_{v}(\tau)e_n$, for $n\in\{0,\dots,2J\}$ and $\tau\geq 0$, the condition $(\rho_{v}(\tau))_{n,n}=0$ is equivalent to $q^n({\tau})=0$. In particular, by assumption, $q^k(\tau)=0$ for $\tau\in[t, t+\varepsilon].$  On this interval we have 
\begin{equation*}
\dot{q}^k(\tau)=iu(\rho_{v}(\tau))\rho_{v}(\tau)J_ye_k=u(\rho_{v}(\tau))\rho_{v}(\tau)\psi = 0,
\end{equation*}
where $\psi := c_k e_{k-1}-c_{k+1} e_{k+1}$. By taking $\varepsilon$ small enough we may assume $u(\rho_{v}(\tau))\neq 0$ and therefore  the previous equality implies $\rho_{v}(\tau)\psi=0.$ This means that $\psi \in Z_2(\tau)$ and, since $\psi\notin Z_1(\tau),$ by the above argument we have $J_z\psi\notin Z_2(\tau)$ and
\begin{equation*}
\psi^*\dot{\rho}_{v}(\tau)\psi=M(1-\eta)\psi^*J_z\rho_{v}(\tau)J_z\psi>0,
\end{equation*}
leading to a contradiction. Hence, there exists $s-t>0$ arbitrarily small such that $Z_1(s)\subsetneqq Z_1(t)$ and, by continuity of $u,$ $u(\rho_{v}(s))\neq 0$. Thus, by repeating the arguments for a finite number of steps, we can show that there exists $s-t>0$ arbitrary small  such that $Z_1(s)=0.$ As $t$ may also be chosen arbitrarily small, this means that there exists an arbitrarily small $s>0$ such that $\rho_{v}(s)>0.$

To conclude the proof, we show that if $\rho_{v}(t_0)>0$ for some $t_0\geq 0,$ then $\rho_{v}(t)>0$ for all $t>t_0.$ This can be done by considering the flow $\Phi_{t,v}: \mathcal{S}\rightarrow\mathcal{S}$ of Equation~\eqref{ODE} which associates with each $\rho_0,$ the value $\rho_{v}(t).$ Since $\Phi_{t,v}$ is a diffeomorphism,  if $\rho\in\mathcal{S}\setminus\partial\mathcal{S}$, there is an open neighborhood $U$ of the state $\rho$ such that $\Phi_{t,v}U \subset \mathcal{S}$ is also an open neighborhood of $\Phi_{t,v}\rho$. Thus, $\Phi_{t,v}\rho\in\mathcal{S}\setminus\partial\mathcal{S}$. The proof is then complete.\hfill$\square$ 
\begin{corollary}
Suppose that the assumptions of Proposition~\ref{Exits boundary ODE} are satisfied. Then for all $\rho_{0}\in\partial\mathcal{S}\setminus{\boldsymbol \rho_{\bar n}}$, either $\rho_t$ stays on the boundary of $\partial \mathcal{S}$ and converges to $\boldsymbol\rho_{\bar n}$ as $t$ goes to infinity or it exits the boundary in finite time  and stays in the interior of $\mathcal{S}$ afterwards, almost surely.  
\label{Exits boundary}
\end{corollary}

\proof
By the support theorem, Theorem~\ref{Support thm}, and Proposition~\ref{Exits boundary ODE}, we have  $\mathbb{P}(\rho_{\nu}>0)>0$ for all $\nu>0$ independently of the initial state $\rho_0\in\mathcal{S}\setminus{\boldsymbol \rho_{\bar n}}.$ Define the set ${\mathcal S}_{\leq \zeta}:=\{\rho\in\mathcal S|\,\det(\rho)\leq\zeta\}\setminus B_r(\boldsymbol\rho_{\bar n})$ for any $r$ arbitrary small and the stopping time $\tau_{\zeta}:=\inf\{t\geq0|\,\rho_t\notin{\mathcal S}_{\leq \zeta} \}$. Now by compactness of ${\mathcal S}_{\leq \zeta}$ and the Feller continuity of $\rho_t$ (\cite[Lemma 4.5]{mirrahimi2007stabilizing}), it is easy to see that for any $\nu>0$ and $\zeta>0$ small enough, there exists $\varepsilon>0$ such that $\mathbb{P}_{\rho_0}(\tau_{\zeta}<\nu)>\varepsilon$,\footnote{Recall that $\mathbb{P}_{\rho_0}$ corresponds to the probability law of $\rho_t$ starting at $\rho_0;$ the associated expectation is denoted by $\mathbb{E}_{\rho_0}$.} independently of $\rho_0\in {\mathcal S}_{\leq \zeta}$. Then we can conclude that $\sup_{\rho_0 \in\mathcal{S}_{\leq\zeta}} \mathbb{P}_{\rho_0}(\tau_{\zeta}\geq \nu) \leq 1-\varepsilon.$
By Dynkin inequality~\cite{dynkin1965markov}, 
\begin{equation*}
\sup_{\rho_0\in{\mathcal S_{\leq \zeta}}}\mathbb{E}_{\rho_0}(\tau_{\zeta}) \leq \frac{\nu}{1-\sup_{\rho_0\in{\mathcal S}_{\leq \zeta}} \mathbb{P}_{\rho_0}(\tau_{\zeta}\geq \nu)}\leq \frac{\nu}{\varepsilon}<\infty.
\end{equation*}
By Markov inequality, for all $\rho_0 \in \mathcal{S}_{\leq\zeta}$, we have 
$$
\mathbb{P}_{\rho_0}(\tau_{\zeta}=\infty) = \lim_{n\rightarrow \infty} \mathbb{P}_{\rho_0}(\tau_{\zeta} \geq n) \leq \lim_{n\rightarrow \infty} \mathbb{E}_{\rho_0}(\tau_{\zeta})/n=0.
$$
By arbitrariness of $r$ we deduce that, either $\rho_t>0$ for some positive time $t$ or $\rho_t$ converges to $\boldsymbol\rho_{\bar n}$ as $t$ tends to infinity while staying in $\partial \mathcal{S}$, almost surely.
In addition, by the strong Markov property of $\rho_t$ and Lemma~\ref{Pos def invariant}, once $\rho_t$ exits the boundary and enters the interior of $\mathcal{S}$, it stays in the interior afterwards. The proof is hence complete.\hfill$\square$

\section{Quantum State Reduction}
\label{sec-qsr}
In this section, we study the dynamics of the $N$-level quantum angular momentum system~\eqref{ND SME} with the feedback $u\equiv0$. First, we can easily show, by Cauchy-Schwarz inequality, that in this case the equilibria of system~\eqref{ND SME} are exactly the eigenstates $\boldsymbol\rho_{n},$ i.e., $F(\boldsymbol\rho_{n})=G(\boldsymbol\rho_{n})=0$ with $n\in\{0,\dots,2J\}$. 

The following theorem shows that the quantum state reduction for the system~\eqref{ND SME} towards the invariant set $\bar E:=\{\boldsymbol\rho_{0},\dots,\boldsymbol\rho_{2J}\}$ occurs with exponential velocity. Note that the exponential stability in mean has been proved independently in the recent paper~\cite{cardona2018exponential}.
  \begin{theorem}[$N$-level quantum state reduction]
For system~\eqref{ND SME}, with $u\equiv0$ and $\rho_0 \in \mathcal{S},$ the set $\bar E$ is exponentially stable in mean and a.s. with average and sample Lyapunov exponent less or equal than $-\eta M/2$. Moreover, the probability of convergence to $\boldsymbol\rho_{n} \in \bar E$ is $\mathrm{Tr}(\rho_0 \boldsymbol\rho_{n})$ for $n \in \{0,\dots,2J\}$.
\label{ND QSR}
\end{theorem}
\proof
Let $I:=\{k| \,\rho_{k,k}(0)=0 \}$ and $\mathcal  S_I:=\{\rho\in \mathcal S| \,\rho_{k,k}=0\mbox{ if and only if } k\in I \}.$ Then by Lemma~\ref{Never reach without feedback}, $\mathcal  S_I$ is a.s. invariant for~\eqref{ND SME}. Consider the function
\begin{equation}
V(\rho)=\frac12\sum^{2J}_{\substack{n,m=0\\n \neq m}}\sqrt{\mathrm{Tr}(\rho\boldsymbol\rho_{n})\mathrm{Tr}(\rho\boldsymbol\rho_{m})} = \frac12\sum^{2J}_{\substack{n,m=0\\n \neq m}}\sqrt{\rho_{n,n}\rho_{m,m}} \geq 0
\label{Lya QSR}
\end{equation}
as a candidate Lyapunov function. Note that $V(\rho)=0$ if and only if $\rho\in\bar{E}$. As $S_I$ is invariant for~\eqref{ND SME} with $u\equiv 0$ and $V$ is twice continuously differentiable when restricted to $S_I$, we can compute
$
\mathscr{L}V(\rho) \leq -\frac{\eta M}{2}V(\rho).
$
By It\^o's formula, for all $\rho_0 \in \mathcal{S}$, we have 
\begin{equation*}
\mathbb{E}(V(\rho_t)) = V(\rho_0)+\int^t_0 \mathbb{E}(\mathscr{L}V(\rho_s))ds \leq V(\rho_0)-\frac{\eta M}{2}\int^t_0 \mathbb{E}(V(\rho_s))ds.
\end{equation*}
In virtue of Gr\"onwall inequality, we have
$
\mathbb{E}(V(\rho_t))\leq V(\rho_0) e^{-\frac{\eta M}{2}t}.
$
Next, we show that the candidate Lyapunov function is bounded by the Bures distance from $\bar E$. Firstly, we have 
\begin{equation*}
V(\rho)=\frac12\sum^{2J}_{n=0}\left(\sqrt{\rho_{n,n}}\sum_{m\neq n}\sqrt{\rho_{m,m}}\right) \geq \frac12\sum^{2J}_{n=0}\sqrt{\rho_{n,n}(1-\rho_{n,n})} \geq \frac{d_B(\rho,\bar{E})}{2}\sum^{2J}_{n=0}\sqrt{\rho_{n,n}}.
\end{equation*}
Combining with $\sum^{2J}_{n=0}\sqrt{\rho_{n,n}} \geq \sum^{2J}_{n=0}\rho_{n,n} =1$, we have 
$
\frac12d_B(\rho,\bar{E}) \leq V(\rho).
$
Let us now prove the converse inequality. Assume that $d_B(\rho,\bar{E}) = \sqrt{2-2\sqrt{\rho_{\bar n,\bar n}}}$  for some index $\bar n,$ then $\sqrt{\rho_{m,m}}\leq \sqrt{1-\rho_{\bar n,\bar n}}\leq  d_B(\rho,\bar{E})$ for $m\neq\bar n$. In particular each addend in $V(\rho)$ is less or equal than $d_B(\rho,\bar{E})$, and $V(\rho)\leq  J(2J+1)d_B(\rho,\bar{E}).$ 

Thus, we have
\begin{equation}
C_1d_B(\rho,\bar{E}) \leq V(\rho) \leq C_2d_B(\rho,\bar{E}),
\label{C1d<=V<=C2d}
\end{equation}
where $C_1 = 1/2$, $C_2= J(2J+1)$. It implies,
\begin{equation*}
\mathbb{E}(d_B(\rho_t,\bar{E})) \leq \frac{C_2}{C_1}d_B(\rho_0,\bar{E})e^{-\frac{\eta M}{2}t}, \quad \forall \rho_0 \in \mathcal{S}.
\end{equation*}
which means that the set $\bar{E}$ is exponentially stable in mean with average Lyapunov exponent less or equal than $-\eta M/2$.

Now we consider the stochastic process
$
Q(\rho_t,t) = e^{\frac{\eta M}{2}t}V(\rho_t) \geq 0
$
whose infinitesimal generator is given by
$
\mathscr{L}Q(\rho,t) = e^{\frac{\eta M}{2}t}( \eta M/2\,V(\rho)+\mathscr{L}V(\rho) )\leq 0.
$ 
Hence, the process $Q(\rho_t,t)$ is a positive supermartingale. Due to Doob's martingale convergence theorem \cite{revuz2013continuous}, the process $Q(\rho_t,t)$ converges almost surely to a finite limit as $t$ tends to  infinity. Consequently, $Q(\rho_t,t)$ is almost surely bounded, that is $\sup_{t \geq 0}Q(\rho_t,t) = A$, for some a.s. finite random variable $A$. This implies $\sup_{t \geq 0} V(\rho_t) = Ae^{-\frac{\eta M}{2}t}$ a.s. Letting $t$ goes  to  infinity, we obtain $\limsup_{t \rightarrow \infty} \frac{1}{t} \log V(\rho_t) \leq -\frac{\eta M}{2}$ a.s. By the inequality~\eqref{C1d<=V<=C2d},
\begin{equation} 
\setlength{\abovedisplayskip}{3pt}
\setlength{\belowdisplayskip}{3pt}
\limsup_{t\rightarrow\infty}\frac{1}{t}\log d_B(\rho_t,\bar E) \leq -\frac{\eta M}{2}, \qquad a.s.
\label{rate QSR}
\end{equation}
which means that the set $\bar E$ is a.s. exponentially stable with sample Lyapunov exponent less or equal than $-\eta M/2$. 

In order to calculate the probability of convergence towards $\boldsymbol\rho_{n} \in \bar E,$ we follow an approach inspired by \cite{amini2013feedback, adler2001martingale}. According to the first part of the theorem, the process
$\mathrm{Tr}(\rho_t\boldsymbol\rho_n)$ converges a.s. to $\mathds{1}_{\{\rho_t\rightarrow\boldsymbol\rho_n\}}.$ Therefore, by applying the dominated convergence theorem, $\mathrm{Tr}(\rho_t\boldsymbol\rho_n)$ converges to $\mathds{1}_{\{\rho_t\rightarrow\boldsymbol\rho_n\}}$ in mean. As $\mathscr{L}\mathrm{Tr}(\rho_t \boldsymbol\rho_{n})=0$, then $\mathrm{Tr}(\rho_t \boldsymbol\rho_{n})$ is a positive martingale. Hence, 
$$
\mathbb{P}(\rho_t\rightarrow\boldsymbol\rho_n)=\lim_{t\rightarrow\infty}\mathbb{E}(\mathrm{Tr}(\rho_t\boldsymbol\rho_{n}))=\mathrm{Tr}(\rho_0\boldsymbol\rho_{n}),
$$
and the proof is complete.\hfill$\square$

\section{Exponential stabilization by continuous feedback}
\label{sec-stab}
In this section, we study the exponential stabilization of system~\eqref{ND SME} towards a selected target state $\boldsymbol\rho_{\bar n}$ with $\bar n\in\{0,\dots,2J\}$. Firstly, we establish a general result ensuring the exponential convergence towards $\boldsymbol\rho_{\bar n}$ under some assumptions on the feedback control law and an additional local Lyapunov type condition. Next, we design  a parametrized family of feedback control laws satisfying such conditions. 

\subsection{Almost sure global exponential stabilization}	  
Inspired by \cite[Lemma 3.4]{tsumura2008global} and \cite[Proposition~3.1]{revuz2013continuous}, in the following lemma we show that, wherever the initial state is, the trajectory $\rho_t$ enters in $B_r(\boldsymbol\rho_{\bar n})$ with $r>0$ in finite time almost surely. 

Before stating the result, we define $\mathbf{P}_{\bar n}:=\{\rho\in\mathcal{S}|\,J-\bar n-\mathrm{Tr}(J_z\rho)=0\}$ and the ``variance function''  $\mathscr V(\rho):=\mathrm{Tr}(J^2_z\rho)-\mathrm{Tr}^2(J_z\rho)$ of $J_z$. 
\begin{lemma}
Assume that $u\in\mathcal{C}^1(\mathcal{S}\setminus\boldsymbol \rho_{\bar n},\mathbb{R})$. Suppose that  for any $\rho_0\in\{\rho\in\mathcal S| \,\, \rho_{\bar n,\bar n}=0\},$  there exists a control $v(t)\in\mathcal{V}$ such that for all $t\in(0,\varepsilon),$ with $\varepsilon$ sufficiently small, $u(\rho_{v}(t))\neq 0$, for some solution $\rho_{v}(t)$ of Equation~\eqref{ODE}. Assume moreover that 
\begin{equation}
\forall \rho\in\mathbf{P}_{\bar n}\setminus \boldsymbol\rho_{\bar n},\quad 2\eta M\mathscr V(\rho)\rho_{\bar n,\bar n}>u(\rho)\mathrm{Tr}(i[J_y,\rho]\boldsymbol\rho_{\bar n}).
\label{Condition u_t}
\end{equation}
Then for all $r>0$ and any given initial state $\rho_0 \in \mathcal{S},$
$
\mathbb{P}(\tau_{r} < \infty)=1,
$
where $\tau_{r}: = \inf\{t \geq 0|\, \rho_t \in B_r(\boldsymbol\rho_{\bar n}) \}$ and $\rho_t$ corresponds to the solution of system~\eqref{ND SME}.
\label{Passage lemma}
\end{lemma}
\proof
The lemma holds trivially for $\rho_{0} \in B_r(\boldsymbol\rho_{\bar n})$, as in that case $\tau_{r} = 0$. Let us thus suppose that $\rho_{0} \in \mathcal{S} \setminus B_r(\boldsymbol\rho_{\bar n})$. We show that there exists $T\in (0,\infty)$ and $\zeta\in (0,1)$ such that $\mathbb{P}_{\rho_0}( \tau_{r}<T )>\zeta$. For this purpose, we make use of the support theorem. Therefore, we consider the  differential equation
\begin{equation}
(\dot{\rho}_{v}(t))_{\bar n,\bar n}=\Delta_{\bar n}(\rho_{v}(t))+2\sqrt{\eta M}P_{\bar{n}}(\rho_{v}(t))(\rho_{v}(t))_{\bar{n},\bar{n}}v(t),
\label{deterministic rho_nn(t)}
\end{equation} 
where $v(t)\in\mathcal{V}$ is the control input,  and
\begin{align*}
\Delta_{\bar n}(\rho)&:=2\eta M\left[\mathrm{Tr}(J^2_z\rho)-(J-\bar{n})^2 \right]\rho_{\bar{n},\bar{n}}-u(\rho)\mathrm{Tr}(i[J_y,\rho]\boldsymbol\rho_{\bar n})\\
&\quad\;+4\eta MP_{\bar{n}}(\rho)\mathrm{Tr}(J_z\rho)\rho_{\bar{n},\bar{n}},\\
P_{\bar{n}}(\rho)&:=J-\bar{n}-\mathrm{Tr}(J_z\rho).
\end{align*} 
Consider the special case in which $\rho_{\bar n,\bar n}(0)=0$. By applying similar arguments as in the proof of Proposition~\ref{Exits boundary ODE}, there exists a control input $v\in\mathcal{V}$ such that $(\rho_{v}(t))_{\bar{n},\bar{n}}>0$ for all $t>0$. Thus, without loss the generality, we suppose $\rho_{\bar n,\bar n}(0)>0$. Then we show that there exist a control input $v$ and a time $T\in(0,\infty)$ such that $\rho_{v}(t)\in B_r(\boldsymbol\rho_{\bar n})$ for $t\leq T$ in the two following separate cases.
\begin{enumerate}
\item Let $\bar{n}\in\{0,2J\}$. We have $\mathbf{P}_{\bar{n}}= \boldsymbol\rho_{\bar n}$. Since $\mathcal{S}\setminus B_r(\boldsymbol \rho_{\bar n})$ is compact, $\Delta_{\bar n}(\rho)$ is bounded from above in this domain and $|{P}_{\bar{n}}(\rho)|$ is bounded from below. 
Then by choosing the control input $v=KP_{\bar{n}}(\rho)/\rho_{\bar n,\bar n}$, with $K>0$ sufficiently large, we can guarantee that $\rho_{v}(t)\in B_r(\boldsymbol\rho_{\bar n})$ for $t\leq T$ with $T<\infty$ if $\rho_{\bar n,\bar n}(0)>0$.
\item Now suppose $\bar{n}\in\{1,\cdots,2J-1\}$. Due to the compactness of $\mathbf{P}_{\bar{n}}\setminus B_r(\boldsymbol \rho_{\bar n})$ and the condition~\eqref{Condition u_t}, we have
\begin{align*}
m:&=\min_{\rho\in\mathbf{P}_{\bar{n}}\setminus B_r(\boldsymbol \rho_{\bar n})}\Delta_{\bar n}(\rho)\\
&=\min_{\rho\in\mathbf{P}_{\bar{n}}\setminus B_r(\boldsymbol \rho_{\bar n})}\Big(2\eta M\mathscr V(\rho)\rho_{\bar n,\bar n}-u(\rho)\mathrm{Tr}(i[J_y,\rho]\boldsymbol\rho_{\bar n})\Big)>0.
\end{align*}
Then we define an open set containing $\mathbf{P}_{\bar{n}}\setminus B_r(\boldsymbol \rho_{\bar n})$,
\begin{equation*}
\mathbf{P}_{\bar{n}}\setminus B_r(\boldsymbol \rho_{\bar n})\subseteq\mathbf{U}:=\{\rho\in\mathcal{S}|\,\Delta_{\bar n}(\rho)>m/2\}\subseteq\mathcal{S}.
\end{equation*} 
Thus, setting $v(t)=0$ whenever $\rho_{v}(t)\in\mathbf{U}$, we have
\begin{equation*}
(\dot{\rho}_{v}(t))_{\bar n,\bar n}=\Delta_{\bar n}(\rho_{v}(t))>m/2\quad\mbox{on }\mathbf{U}.
\end{equation*}
Moreover, $(\mathcal{S}\setminus B_r(\boldsymbol \rho_{\bar n}))\setminus\mathbf{U}$ is compact, then $\Delta_{\bar n}(\rho)$ is bounded from above and $|P_{\bar n}(\rho)|$ is bounded from below in this domain. For all $\rho_{v}(t) \in \{\rho\in\mathcal S|\,\rho_{\bar{n},\bar{n}}>0\}$, we can take the feedback $v=KP_{\bar{n}}(\rho)/\rho_{\bar n,\bar n}$ with $K>0$ sufficiently large, so that $(\dot{\rho}_{v}(t))_{\bar n,\bar n}$ is bounded from below on $(\mathcal{S}\setminus B_r(\boldsymbol \rho_{\bar n}))\setminus\mathbf{U}$. The proposed input $v$ guarantees that  
 $\rho_{v}(t)\in B_r(\boldsymbol \rho_{\bar n})$ for $t\leq T$ with $T<\infty$ if $\rho_{\bar n,\bar n}(0)>0$.
\end{enumerate}
Therefore, there exists  $T\in (0,\infty)$ such that, for all $\rho_{0}\in \mathcal S\setminus B_r(\boldsymbol\rho_{\bar n})$, there exists  $v(t)$ steering the system  from $\rho_0$ to $B_r(\boldsymbol\rho_{\bar n})$ by time $T.$ By compactness of $\mathcal S\setminus B_r(\boldsymbol\rho_{\bar n})$ and the Feller continuity of $\rho_t,$ we have 
$
\sup_{\rho_0 \in\mathcal S\setminus B_r(\boldsymbol\rho_{\bar n}) } \mathbb{P}_{\rho_0}(\tau_{r}\geq T) \leq 1-\zeta<1,
$
for some $\zeta>0.$
By Dynkin inequality~\cite{dynkin1965markov}, 
\begin{equation*}
\sup_{\rho_0 \in\mathcal S\setminus B_r(\boldsymbol\rho_{\bar n}) } \mathbb{E}_{\rho_0}(\tau_{r}) \leq \frac{T}{1-\sup_{\rho_0 \in\mathcal S\setminus B_r(\boldsymbol\rho_{\bar n}) } \mathbb{P}_{\rho_0}(\tau_{r}\geq T)}\leq \frac{T}{\zeta}<\infty.
\end{equation*}
Then by Markov inequality, for all $\rho_0 \in\mathcal S\setminus B_r(\boldsymbol\rho_{\bar n})$, we have 
\begin{equation*}
\mathbb{P}_{\rho_0}(\tau_{r}=\infty) = \lim_{n\rightarrow \infty} \mathbb{P}_{\rho_0}(\tau_{r} \geq n) \leq \lim_{n\rightarrow \infty} \mathbb{E}_{\rho_0}(\tau_{r})/n=0,
\end{equation*}
which implies
$
\mathbb{P}_{\rho_0}( \tau_{r}<\infty )=1.
$
The proof is  complete.\hfill$\square$

\medskip
In the following, we state our general result concerning the exponential stabilization of $N$-level quantum angular momentum systems. 
\begin{theorem}
Assume that the feedback control law satisfies the assumptions of Lemma~\ref{Never reach with feedback} and Lemma~\ref{Passage lemma}. Additionally, suppose that there exists a positive-definite function $V(\rho)$ such that $V(\rho)=0$ if and only if $\rho=\boldsymbol\rho_{\bar n}$, and $V$ is continuous on $\mathcal{S}$ and twice continuously differentiable on the set $\mathcal S\setminus\boldsymbol\rho_{\bar n}$. Moreover, suppose that there exist positive constants $C$, $C_1$ and $C_2$ such that 
\begin{enumerate}
\item[(i)] $C_1 \, d_B(\rho,\boldsymbol\rho_{\bar n}) \leq V(\rho) \leq C_2 \, d_B(\rho,\boldsymbol\rho_{\bar n})$, for all $\rho\in\mathcal{S}$, and 
\item[(ii)] $\limsup_{\rho\rightarrow\boldsymbol\rho_{\bar n}}\frac{\mathscr{L}V(\rho)}{V(\rho)}=-C$.
\end{enumerate}
Then, $\boldsymbol\rho_{\bar n}$ is a.s. exponentially stable for the system~\eqref{ND SME} with sample Lyapunov exponent less or equal than $-C-\frac{K}{2}$, where $K:=\liminf_{\rho \rightarrow\boldsymbol \rho_{\bar{n}}}g^2(\rho)$ and $g(\rho):=\sqrt{\eta}\frac{\partial V(\rho)}{\partial \rho}\frac{G(\rho)}{V(\rho)}$.
\label{Thm a.s. exp stab 1}
\end{theorem}
\proof
The proof proceeds in three steps:
\begin{enumerate}
\item First we show that $\boldsymbol\rho_{\bar n}$ is locally stable in probability;
\item Next we show that for any fixed $r>0$ and almost all sample paths, there exists $T<\infty$ such that  for all $t\geq T$, $\rho_t \in B_r(\boldsymbol\rho_{\bar n})$;
\item Finally, we prove that $\boldsymbol\rho_{\bar n}$ is a.s. exponentially stable with sample Lyapunov exponent less or equal than $-C-\frac{K}{2}$.
\end{enumerate}
\textbf{Step 1:}
By the condition \textit{(ii)}, we can choose $r>0$ sufficiently small such that $\mathscr{L}V(\rho)\leq -C(r)V(\rho)$ for $\rho\in B_r(\boldsymbol\rho_{\bar n})\setminus\boldsymbol \rho_{\bar n},$ for some $C(r)>0.$ Let $\varepsilon \in (0,1)$ be arbitrary. By the continuity of $V(\rho)$ and the fact that $V(\rho)=0$ if and only if $d_B(\rho,\boldsymbol\rho_{\bar n})=0$, we can find  $\delta=\delta(\varepsilon,r)>0$ such that
\begin{equation}
1/\varepsilon\sup_{\rho_{0} \in B_{\delta}(\boldsymbol\rho_{\bar n})}V(\rho_{0}) \leq C_1r.
\label{Step1:stable ineq}
\end{equation} 
Assume that $\rho_{0} \in B_{\delta}(\boldsymbol\rho_{\bar n})$ and  let $\tau$ be the first exit time of $\rho_t$ from $B_r(\boldsymbol\rho_{\bar n})$. By It\^o's formula, we have 
\begin{equation*}
\mathbb{E}(V(\rho_{t \wedge \tau})) \leq V(\rho_{0})-C(r) \, \mathbb{E} \left(\int^{t \wedge \tau}_{0} V(\rho_s)ds \right) \leq V(\rho_{0}).
\end{equation*}
For all $t \geq \tau$, $d_B(\rho_{t \wedge \tau},\boldsymbol\rho_{\bar n})=d_B(\rho_{\tau},\boldsymbol\rho_{\bar n})=r$. Hence, by the condition~\textit{(i)},
\begin{equation*}
\mathbb{E}(V(\rho_{t \wedge \tau})) \geq \mathbb{E}(\mathds{1}_{\{\tau\leq t\}}V(\rho_{\tau})) \geq \mathbb{E}(\mathds{1}_{\{\tau\leq t\}}C_1d_B(\rho_{\tau},\boldsymbol\rho_{\bar n})) = C_1 r \, \mathbb{P}(\tau \leq t).
\end{equation*}
Combining with the inequality~\eqref{Step1:stable ineq}, we have 
\begin{equation*}
\mathbb{P}(\tau \leq t) \leq \frac{\mathbb{E}(V(\rho_{t \wedge \tau}))}{C_1r} \leq \frac{V(\rho_{0})}{C_1r} \leq \varepsilon.
\end{equation*}
Letting $t$ tend to infinity, we get $\mathbb{P} (\tau<\infty) \leq \varepsilon$ which  implies 
\begin{equation*}
\mathbb{P}(d_B(\rho_t,\boldsymbol\rho_{\bar n}) < r \text{ for } t \geq 0) \geq 1-\varepsilon.
\end{equation*}
\textbf{Step 2:}
Since $u_t=0$ in $\bar E$ if and only if $\rho_t=\boldsymbol\rho_{\bar n}$ by Lemma~\ref{Passage lemma} we obtain, for all $\rho_0 \in \mathcal{S}$, $\mathbb{P}(\tau_{\delta} < \infty)=1$, where $\tau_{{\delta}}: = \inf\{t \geq 0|\, \rho_t \in B_{\delta}(\boldsymbol\rho_{\bar n}) \}$. It implies that $\rho_t$ enters $B_{\delta}(\boldsymbol\rho_{\bar n})$ in a finite time almost surely. Due to Step 1, for all $\rho_0 \in B_{\delta}(\boldsymbol\rho_{\bar n})$, $\mathbb{P}(\sigma_{r} < \infty)\leq \varepsilon$, where $\sigma_{r} := \inf\{t \geq 0|\,\rho_t \notin B_{r}(\boldsymbol\rho_{\bar n})\}$. 

We define two sequences of stopping times $\{\sigma^k_r\}_{k\geq 0}$ and $\{\tau^k_\delta\}_{k\geq 1}$ such that $\sigma^0_r=0,$  $\tau^{k+1}_\delta = \inf\{t \geq \sigma^k_r|\, \rho_t \in B_{\delta}(\boldsymbol\rho_{\bar n})\}$ and $\sigma^{k+1}_r = \inf\{t \geq \tau^{k+1}_\delta|\,\rho_t \notin B_r(\boldsymbol\rho_{\bar n})\}$. By the strong Markov property, we find
\begin{equation*}
\begin{split}
\mathbb{P}_{\rho_0}(\sigma^m_r<\infty)&=\mathbb{P}_{\rho_0}(\tau^1_\delta<\infty,\sigma^1_r<\infty,\dots,\sigma^m_r<\infty)\\
&=\mathbb{P}_{\rho_{\tau^1_\delta}}(\sigma_{r}<\infty)\cdots\mathbb{P}_{\rho_{\tau^m_\delta}}(\sigma_{r}<\infty)
\leq \varepsilon^m.
\end{split}
\end{equation*}
Thus,  for all $\rho_0\in\mathcal S,$ we have
$\mathbb{P}(\sigma^m_r<\infty,\,\forall m>0)=0.
$
We deduce that, for almost all sample paths, there exists $T<\infty$ such that, for all $t \geq T$, $\rho_t \in B_r(\boldsymbol\rho_{\bar n})$, which concludes Step~2. 

\textbf{Step 3:}
In this step, we obtain an upper bound of the sample Lyapunov exponent by employing an argument inspired by~\cite[Theorem 4.3.3]{mao2007stochastic}. For $\rho\neq\boldsymbol \rho_{\bar n}$, 
$
\mathscr{L}\log V(\rho) = \frac{\mathscr{L}V(\rho)}{V(\rho)}-\frac{g^2(\rho)}{2}.
$
Due to Lemma~\ref{Never reach with feedback}, $\boldsymbol \rho_{\bar n}$ cannot be attained in finite time almost surely, then by It\^o's formula, we have
\begin{equation*}
\log V(\rho_t)=\log V(\rho_{0})+ \int^t_0\frac{\mathscr{L}V(\rho_s)}{V(\rho_s)}ds+\int^t_0g(\rho_s)dW_s -\frac12\int^t_0g^2(\rho_s)ds.
\end{equation*}
Let $m\in\mathbb Z_{>0}$ and take arbitrarily $\varepsilon\in(0,1)$. By the exponential martingale inequality (see e.g.~\cite[Theorem 1.7.4]{mao2007stochastic}), we have
\begin{equation*}
\mathbb{P} \left(\sup_{0 \leq t\leq m} \!\left[ \int^t_{0}\!\!g(\rho_s)dW_s-\frac{\varepsilon}{2}\!\int^t_{0}\!\!g^2(\rho_s)ds \right] \!> \!\frac{2}{\varepsilon} \log m \right)\leq\! \frac{1}{m^2}.
\end{equation*}
Since $\sum^{\infty}_{m=1}\frac{1}{m^2}<\infty$, by Borel-Cantelli lemma we have that for almost all sample paths there exists $m_0$ such that, if $m>m_0$, then
\begin{equation*}
\sup_{0\leq t\leq m}\left( \int^t_{0}g(\rho_s)dW_s-\frac{\varepsilon}{2}\int^t_{0}g^2(\rho_s)ds \right)\leq\frac{2}{\varepsilon}\log m.
\end{equation*}
Thus, for $0\leq t\leq m$ and $m>m_0$,
\begin{equation*}
\int^t_{0}g(\rho_s)dW_s \leq \frac{2}{\varepsilon}\log m+\frac{\varepsilon}{2}\int^t_{0}g^2(\rho_s)ds,\quad a.s.
\end{equation*}
We have
\begin{equation*}
\log V(\rho_t) \leq \log V(\rho_0)+\int^t_0\frac{\mathscr{L}V(\rho_s)}{V(\rho_s)}ds+\frac{2}{\varepsilon}\log m-\frac{1-\varepsilon}{2}\int^t_{0}g^2(\rho_s)ds,\quad a.s.
\end{equation*}
It gives
\begin{equation*} 
\limsup_{t \rightarrow \infty}\frac{1}{t}\log V(\rho_t) \leq \limsup_{t \rightarrow \infty}\frac{1}{t}\left(\int^t_0\frac{\mathscr{L}V(\rho_s)}{V(\rho_s)}ds-\frac{1-\varepsilon}{2}\int^t_{0}g^2(\rho_s)ds\right) \quad a.s.
\end{equation*}
Letting $\varepsilon$ tend to zero, we have 
\begin{align*}
\limsup_{t \rightarrow \infty}\frac{1}{t}\log V(\rho_t) \leq \limsup_{t \rightarrow \infty}\frac{1}{t}\left(\int^t_0\frac{\mathscr{L}V(\rho_s)}{V(\rho_s)}ds- \frac{1}{2}\int^t_{0}g^2(\rho_s)ds\right)\quad a.s.
\end{align*}
For every fixed $T>0$ consider the event 
\begin{equation*}
\Omega_T=\{\rho_t \in  B_r(\boldsymbol \rho_{\bar{n}}) \text{ for all } t\geq T\}.
\end{equation*}
Due to the condition \textit{(ii)}, for almost all $\omega\in\Omega_T$, 
\begin{align*}
\limsup_{t \rightarrow \infty}\frac{1}{t}&\left(\int^t_0\frac{\mathscr{L}V(\rho_s)}{V(\rho_s)}ds- \frac{1}{2}\int^t_{0}g^2(\rho_s)ds\right)\\ 
&\qquad\qquad\leq\limsup_{t\rightarrow\infty}\frac{1}{t}\left(\int^t_T\frac{\mathscr{L}V(\rho_s)}{V(\rho_s)}ds- \frac{1}{2}\int^t_{T}g^2(\rho_s)ds\right)\\
&\qquad\qquad\leq -C(r)-\inf_{\rho \in  B_r(\boldsymbol \rho_{\bar{n}})\setminus\boldsymbol \rho_{\bar{n}}}\frac{g^2(\rho)}{2}.
\end{align*}
Since $T$ can be taken arbitrarily large and Step 2 implies that $\lim_{T\to \infty}\mathbb{P}(\Omega_T)=1$, we can conclude that 
\begin{equation*}
\limsup_{t\rightarrow\infty}\frac{1}{t}\log V(\rho_t) \leq -C(r)-\inf_{\rho \in  B_r(\boldsymbol \rho_{\bar{n}})\setminus\boldsymbol \rho_{\bar{n}}}\frac{g^2(\rho)}{2},\quad a.s.
\end{equation*}
Finally, due to the condition \textit{(i)} and since $r$ can be taken arbitrarily small, we have  
\begin{equation*}
\limsup_{t\rightarrow\infty}\frac{1}{t}\log d_B(\rho_t,\boldsymbol \rho_{\bar{n}}) \leq -C-\frac{K}{2},\quad a.s.
\end{equation*}
which yields the result. \hfill$\square$

\subsection{Feedback controller design}
The purpose of this subsection is to design parametrized feedback laws which stabilize exponentially  the system~\eqref{ND SME} almost surely towards some predetermined target eigenstate. For the choice of target state, we consider first the particular case $\bar n\in\{0,2J\}$ and then the general case $\bar n\in\{0,\cdots,2J\}.$ 

In the following theorem, we consider the case $\bar n\in\{0,2J\}.$ Before stating the result, we note that we can describe the set $B_{r(\lambda)}(\boldsymbol \rho_{\bar{n}})\setminus \boldsymbol \rho_{\bar{n}}$ as follows
\begin{equation*}
D_{\lambda}(\boldsymbol \rho_{\bar{n}}) := \{\rho \in \mathcal{S}|\,0 < \lambda < \rho_{\bar{n},\bar{n}} < 1\} = B_{r(\lambda)}(\boldsymbol \rho_{\bar{n}})\setminus \boldsymbol \rho_{\bar{n}},
\end{equation*}
where $r(\lambda) = \sqrt{2-2\sqrt{\lambda}}$.
\begin{theorem}
Consider system~\eqref{ND SME} with $\rho_{0} \in \mathcal{S}$ and assume $\eta\in(0,1)$. Let $\boldsymbol \rho_{\bar{n}} \in \{\boldsymbol \rho_{0},\boldsymbol \rho_{2J}\}$ be the target eigenstate and define the feedback controller
\begin{equation}
u_{\bar{n}}(\rho)  = \alpha (1-\mathrm{Tr}(\rho \boldsymbol \rho_{\bar{n}}))^{\beta} - \gamma \, \mathrm{Tr}(i[J_y,\rho]\boldsymbol \rho_{\bar{n}}),
\label{u_t 0 2J}
\end{equation}
where $\gamma \geq 0$, $\beta > 1/2$ and $\alpha >0 $. Then the feedback controller~\eqref{u_t 0 2J} exponentially stabilizes system~\eqref{ND SME} almost surely to the equilibrium $\boldsymbol \rho_{\bar{n}}$ with sample Lyapunov exponent less or equal than $-\eta M$.
\end{theorem}
\proof
To prove the theorem, we show that we can apply Theorem~\ref{Thm a.s. exp stab 1} with the Lyapunov function 
$
V_{\bar{n}}(\rho) = \sqrt{1-\mathrm{Tr}(\rho \boldsymbol \rho_{\bar{n}})}
$
for $\bar n=0$ and $\bar n=2J.$ First, it is easy to see that $u_{\bar n}$ satisfies the assumptions of Lemma~\ref{Passage lemma} and Lemma~\ref{Never reach with feedback}. Then, we need to show that the conditions \textit{(i)} and \textit{(ii)} of Theorem~\ref{Thm a.s. exp stab 1} hold true. Note that $\frac{\sqrt{2}}{2}d_B(\rho,\boldsymbol \rho_{\bar{n}}) \leq V_{\bar{n}}(\rho) \leq d_B(\rho,\boldsymbol \rho_{\bar{n}}),$ so that the condition  \textit{(i)}  is shown. We are left to check the condition  \textit{(ii)}.
The infinitesimal generator $\mathscr{L}V_{\bar{n}}$ takes the following form
\begin{equation}
\begin{split}
\mathscr{L}V_{\bar{n}}(\rho)=\frac{u_{\bar{n}}}{2}\frac{\mathrm{Tr}(i[J_y,\rho]\boldsymbol \rho_{\bar{n}})}{V_{\bar{n}}(\rho)}-\frac{\eta M}{2}\frac{(J-\bar{n}-\mathrm{Tr}(J_z\rho))^2\mathrm{Tr}^2(\rho\boldsymbol \rho_{\bar{n}})}{V^3_{\bar{n}}(\rho)}.
\end{split}
\label{LV 0 2J}
\end{equation}
If $\bf \bar{n}=0$, and $\rho \in D_{\lambda}(\boldsymbol \rho_0),$ we find 
$$
\frac{u_{0}}{2}\frac{\mathrm{Tr}(i[J_y,\rho]\boldsymbol \rho_0)}{V_0(\rho)} \leq \alpha c_1 (V_0(\rho))^{\beta} \leq \alpha c_1 (1-\lambda)^{\frac{\beta-1}{2}}V_0(\rho),
$$
since
$
|\mathrm{Tr}(i[J_y,\rho]\boldsymbol \rho_0)|=2c_1|\mathbf{Re}\{\rho_{0,1}\}| \leq 2c_1|\rho_{0,1}| \leq 2c_1V_0(\rho).
$
Moreover, we have $$J-\mathrm{Tr}(J_z\rho)=\sum^{2J}_{k=1}k\rho_{k,k}\geq \sum^{2J}_{k=1}\rho_{k,k}=1-\rho_{0,0}=(V_0(\rho))^2.$$
Thus, for all $\rho \in D_{\lambda}(\boldsymbol \rho_0)$, 
$
\mathscr{L}V_0(\rho)\leq -C_{0,\lambda}V_0(\rho), 
$
where $C_{0,\lambda}=\frac{\eta M \lambda^2}{2}-\alpha c_1 (1-\lambda)^{\frac{\beta-1}{2}}$.
The case  $\bf \bar{n}=2J$ may be treated similarly. In particular,  for all $\rho \in D_{\lambda}(\boldsymbol \rho_{2J})$, one gets
$
\mathscr{L}V_{2J}(\rho)\leq -C_{2J,\lambda}V_{2J}(\rho),
$
where $C_{2J,\lambda}=\frac{\eta M \lambda^2}{2}-\alpha c_{2J} (1-\lambda)^{\frac{\beta-1}{2}}=C_{0,\lambda}$.

Furthermore, for $\bar{n} \in \{0,2J\},$ we have 
$
g^2(\rho)  \geq \eta M \lambda^2,
$ for all $\rho \in D_{\lambda}(\boldsymbol \rho_{\bar{n}}).$
Hence, we can apply Theorem~\ref{Thm a.s. exp stab 1} for $\bar n\in\{0, 2J\},$ with $C=\frac{\eta M}{2}$ and $K = \eta M.$  
The proof is  complete.\hfill$\square$
\medskip

In the following theorem, we consider the general case $\bar n\in\{0,\dots,2J\}.$
\begin{theorem}
Consider system~\eqref{ND SME} with $\rho_{0} \in \mathcal{S}\setminus\partial \mathcal{S}$.  Let $\boldsymbol \rho_{\bar{n}} \in \bar{E}$ be the target eigenstate and define the feedback 
\begin{equation}
u_{\bar{n}}(\rho) = \alpha (P_{\bar{n}}(\rho))^\beta = \alpha (J-\bar{n}-\mathrm{Tr}(J_z\rho))^{\beta},
\label{u_t N}
\end{equation}
where $\beta>1/2$ and $\alpha>0$. Then the feedback~\eqref{u_t N} exponentially stabilizes system~\eqref{ND SME} almost surely to the equilibrium $\boldsymbol\rho_{\bar{n}}$ with  sample Lyapunov exponent less or equal than $-\eta M$ if $\bar{n}\in\{0,2J\}$ and $-\eta M/2$ if $\bar{n}\in\{1,\dots,2J-1\}$.
\label{Th feedback N}
\end{theorem}
\proof
Consider the following candidate Lyapunov function
\begin{equation}
V_{\bar{n}}(\rho) = \sum_{k\neq \bar{n}}\sqrt{\mathrm{Tr}(\rho\boldsymbol \rho_k)}.
\label{Lya fun N}
\end{equation}
Due to Lemma~\ref{Pos def invariant}, all diagonal elements of $\rho_t$ remain strictly positive for all $t\geq0$ almost surely. Since $V_{\bar{n}}(\rho)$ is $\mathcal C^2$ in $\mathcal{S}\setminus\partial \mathcal{S},$ we can make use of similar arguments as those in Theorem~\ref{Thm a.s. exp stab 1}. 
First, we show that the  following conditions are satisfied. 
\begin{enumerate}
\item[C.1.] $2\eta M \mathscr V(\rho)\rho_{\bar{n},\bar{n}}>u_{\bar n}\mathrm{Tr}(i[J_y,\rho]\boldsymbol \rho_{\bar{n}})$, $\forall\,\rho\in\mathbf{P}_{\bar n}\setminus\boldsymbol \rho_{\bar{n}}$,
\item[C.2.] $\sum_k\frac{\partial P_{\bar{n}}(\rho)}{\partial \rho_{k,k}}(G(\rho))_{k,k}\neq 0$ when $u_{\bar{n}}(\rho)=0$ and $\rho\neq\boldsymbol\rho_{\bar{n}}$,
\item[C.3.] $u_{\bar n}(\rho)\leq\overline{C}V_{\bar{n}}(\rho)$ with $\overline{C}>0$, $\forall\,\rho\in D_{\lambda}(\boldsymbol \rho_{\bar{n}})$.
\end{enumerate}  
Roughly speaking, C.1 and C.2 ensure that the assumptions of Proposition~\ref{Exits boundary ODE}, and so those of Lemma~\ref{Passage lemma}, hold true; in particular, C.1 provides a sufficient condition guaranteeing the accessibility  of any arbitrary small neighborhood of $\boldsymbol \rho_{\bar{n}}$. C.3 is helpful to obtain a bound of the type $\mathscr{L}V_{\bar{n}}\leq -CV_{\bar{n}}$ on $D_{\lambda}(\boldsymbol \rho_{\bar{n}})$.

We now show that these conditions are satisfied. The property C.1 follows from the fact that,  for all $\rho \in \mathbf{P}_{\bar{n}}\setminus\boldsymbol \rho_{\bar{n}}$, we have $u_{\bar n}(\rho)=0$ and $\mathscr V(\rho)>0$. 

The condition C.2 can be proved by  contradiction as follows. We suppose  $u_{\bar{n}}(\rho)=0$, $\rho\neq\boldsymbol\rho_{\bar{n}}$ and $\sum_k\frac{\partial P_{\bar{n}}(\rho)}{\partial \rho_{k,k}}(G(\rho))_{k,k}=0.$ Then  it is not difficult to see that $\mathrm{Tr}(J^2_z\rho)=(J-\bar{n})^2=(\mathrm{Tr}(J_z\rho))^2$, that is $\mathscr V(\rho)=0.$ By applying Cauchy-Schwarz inequality, this implies that $\rho\in \bar{E}\setminus\boldsymbol\rho_{\bar{n}},$ which contradicts the fact that $u_{\bar{n}}(\rho)=0.$ 

Finally, we can show that the property C.3 holds true, because
\begin{equation*}
|P_{\bar{n}}(\rho)|=\Big|\sum_{k\neq \bar{n}}k\rho_{k,k}-\bar{n}(1-\rho_{\bar{n},\bar{n}})\Big|\leq \Upsilon(1-\rho_{\bar{n},\bar{n}}),
\end{equation*}
where $\Upsilon := \max\{\bar{n},2J-\bar{n}\}$. Then, for all $\rho\in D_{\lambda}(\boldsymbol \rho_{\bar{n}})$, 
\begin{equation*}
u_{\bar n}(\rho)\leq \alpha\Upsilon^{\beta}(1-\rho_{\bar{n},\bar{n}})^{\beta-1/2}\sqrt{1-\rho_{\bar{n},\bar{n}}}\leq \alpha\Upsilon^{\beta}(1-\lambda)^{\beta-1/2}V_{\bar{n}}(\rho).
\end{equation*} 

Consider the Lyapunov function~\eqref{Lya fun N}. In the following, we verify the conditions \textit{(i)} and \textit{(ii)} of Theorem~\ref{Thm a.s. exp stab 1}. First note that by Jensen's inequality, we have $V_{\bar{n}}(\rho) \leq \sqrt{2J}\sqrt{1-\rho_{\bar{n},\bar{n}}}.$ Then we get $\frac{\sqrt{2}}{2}d_B(\rho,\boldsymbol \rho_{\bar{n}}) \leq V_{\bar{n}}(\rho) \leq \sqrt{2J}d_B(\rho,\boldsymbol \rho_{\bar{n}}),$ hence the condition \textit{(i)} is shown. In order to verify the condition \textit{(ii)}, we write the  infinitesimal generator of the Lyapunov function which has the following form
\begin{equation}
\begin{split}
\mathscr{L}V_{\bar{n}}(\rho) = -\frac{u_{\bar{n}}}{2}\sum_{k\neq \bar{n}}\frac{\mathrm{Tr}(i[J_y,\rho]\boldsymbol \rho_{k})}{\sqrt{\rho_{k,k}}}-\frac{\eta M}{2} \sum_{k\neq \bar{n}} (P_k(\rho))^2 \sqrt{\rho_{k,k}}.
\end{split}
\label{LV N}
\end{equation}
We find
\begin{equation*}
\begin{split}
\frac{|\mathrm{Tr}(i[J_y,\rho]\boldsymbol \rho_{k})|}{\sqrt{\rho_{k,k}}} &= \frac{|c_k\mathbf{Re}\{\rho_{k,k-1}\}-c_{k+1}\mathbf{Re}\{\rho_{k,k+1}\}|}{\sqrt{\rho_{k,k}}}\leq \frac{c_k|\rho_{k,k-1}|+c_{k+1}|\rho_{k,k+1}|}{\sqrt{\rho_{k,k}}}\\
& \leq  c_k\sqrt{\rho_{k-1,k-1}}+c_{k+1}\sqrt{\rho_{k+1,k+1}} \leq c_k+c_{k+1}.
\end{split}
\end{equation*}
 For $k\neq\bar n$ and for all $\rho \in D_{\lambda}(\boldsymbol \rho_{\bar{n}})$ with $\lambda > 1-1/\Upsilon$, we have
\begin{equation*}
|J-k-\mathrm{Tr}(J_z\rho)|\geq  |\bar{n}-k|-|P_{\bar{n}}(\rho)| \geq 1-\Upsilon (1-\rho_{\bar{n},\bar{n}}) \geq 1-\Upsilon(1-\lambda)>0.
\end{equation*}
Thus, for all $\rho \in D_{\lambda}(\boldsymbol \rho_{\bar{n}})$, 
\begin{equation*}
\mathscr{L}V_{\bar{n}}(\rho) \leq -\left( \frac{\eta M(1-\Upsilon(1-\lambda))^2}{2}-\alpha\Gamma\Upsilon^{\beta}(1-\lambda)^{\beta-1/2} \right) V_{\bar{n}}(\rho) \leq -C_{\bar{n},\lambda}V_{\bar{n}}(\rho),
\end{equation*}
where $\Gamma:=\sum_{k\neq \bar{n}}(c_k+c_{k+1})$ and $C_{\bar{n},\lambda} := \frac{\eta M(1-\Upsilon(1-\lambda))^2}{2}-\alpha\Gamma\Upsilon^{\beta}(1-\lambda)^{\beta-1/2}.$

Furthermore, for $\bar{n}\in\{0,2J\}$, we have
$
g^2(\rho) \geq \eta M\lambda^2,
$ for all $\rho \in D_{\lambda}(\boldsymbol \rho_{\bar{n}}).$ Since $C_{\bar{n},\lambda}$ and $\eta M\lambda^2$ converge respectively to $\frac{\eta M}{2}$ and $\eta M$ as $\lambda$ tends to one, by employing the same arguments used earlier in the proof of Theorem~\ref{Thm a.s. exp stab 1}, we find that the sample Lyapunov exponent is less or equal than $-C-K/2$ where  $C=\frac{\eta M}{2}$ for $\bar{n}\in\{0,\dots,2J\},$ $K = \eta M$ for  $\bar{n}\in\{0,2J\}$ and $K=0$ for $\bar{n}\in\{1,\dots,2J-1\}$. \hfill$\square$
\begin{remark}
Locally around the target eigenstate $\boldsymbol \rho_{\bar{n}}$, the asymptotic behavior of the Lyapunov function~\eqref{Lya fun N} is the same as the one of the Lyapunov function~\eqref{Lya QSR}. This is related to the fact that, under the assumptions on $u_{\bar n}$, the behavior of the system around the target state is similar to the case $u\equiv 0.$  In particular, without feedback and conditioning to the event $\{\exists t'\geq 0|\,\rho_t \in B_r(\boldsymbol \rho_{\bar{n}}),\;\forall t\geq t'\}$, one can show that the trajectories converge a.s. to $\boldsymbol \rho_{\bar{n}}$ with sample Lyapunov exponent equal to the one in Theorem~\ref{Th feedback N}.
\end{remark}

\begin{remark}
If $\eta\in(0,1),$  Theorem~\ref{Th feedback N} and Corollary~\ref{Exits boundary} guarantee the convergence of almost all trajectories to the target state even if the initial state $\rho_0$ lies in the boundary of $\mathcal S$ (the argument is no more valid if $\eta=1$ because of Lemma~\ref{boundary}). Unfortunately, these results do not ensure the almost sure exponential convergence towards the target state whenever $\rho_0$ lies in $\partial \mathcal S\setminus\boldsymbol\rho_{\bar n}.$ However, we believe that under the assumptions imposed on the feedback, we can still guarantee such convergence property. This is suggested by the following arguments. 
 
Set the event $\Omega_{>0}=\bigcap_{t>0}\{\rho_t>0\}$ which is $\mathcal{F}_{0+}$-measurable. By the strong Markov property of $\rho_t$, and by applying Blumenthal's zero--one law~\cite{rogers2000diffusions1}, we have that either $\mathbb{P}(\Omega_{>0})=0$ or $\mathbb{P}(\Omega_{>0})=1$. In order to conclude that $\mathbb{P}(\Omega_{>0})=1,$ it would be enough to show that $\mathbb{P}(\Omega_{>0})>0$, i.e., $\rho_t$ exits the boundary and enters the interior of $\mathcal{S}$ immediately with non-zero probability. Proposition~\ref{Exits boundary ODE} provides some intuitions about the validity of this property, as it proves that the majority of the trajectories of the associated deterministic equation~\eqref{ODE} enter the interior of $\mathcal{S}$ immediately. It is then tempting to conjecture that under the assumption of Proposition~\ref{Exits boundary ODE}, for all $\rho_0\in\partial\mathcal{S}\setminus\boldsymbol \rho_{\bar n}$, $\rho_t>0$ for all $t>0$ almost surely. If this conjecture is correct, we can generalize Theorem~\ref{Th feedback N} to the case $\rho_0\in\mathcal{S}.$
\end{remark}

\section{Simulations}
In this section, we illustrate our results by numerical simulations in the case of a three-level quantum angular momentum system. First, we consider the case $u\equiv 0$ (Theorem~\ref{ND QSR}). Then, we illustrate the convergence towards the target states $\boldsymbol\rho_{0}$ and $\boldsymbol\rho_{1}$ by applying feedback laws of the form~\eqref{u_t 0 2J} and~\eqref{u_t N}, respectively. 

\begin{figure}[thpb]
\centering
      \includegraphics[width=12cm]{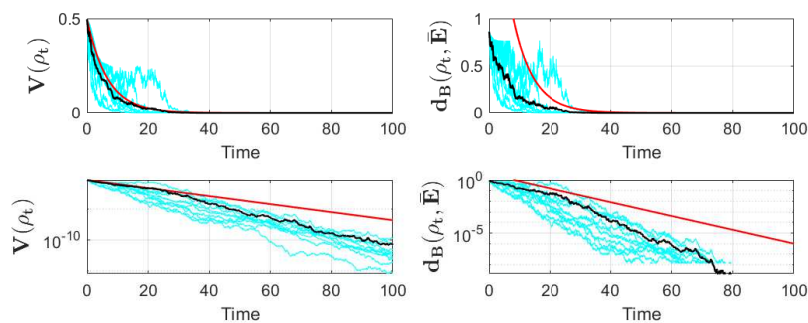}
      \caption{Quantum state reduction of a three-level quantum angular momentum system with $u\equiv 0$ starting at $diag(0.3,0.4,0.3)$ when $\omega=0$, $\eta=0,3$ and $M=1$: the black curve represents the mean value of  10 arbitrary sample trajectories, the red curve represents the exponential reference with exponent $-\eta M/2$. The figures at the bottom are the semi-log versions of the ones at the top. }
      \label{QSR3fig}
\end{figure}

\begin{figure}[thpb]
\centering
      \includegraphics[width=12.5cm]{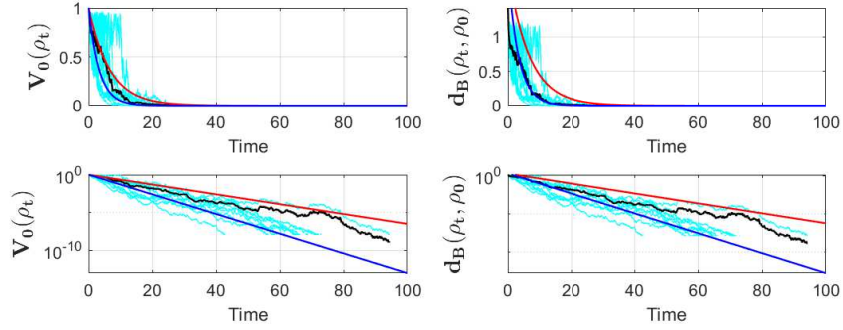}
      \caption{Exponential stabilization of a three-level quantum angular momentum system towards $\boldsymbol \rho_0$ with the feedback law~\eqref{u_t 0 2J}starting at $\boldsymbol \rho_2$ with $\omega=0$, $\eta=0.3$, $M=1$, $\alpha = 10$, $\beta = 5$ and $\gamma = 10$: the black curve represents the mean value of  10 arbitrary sample trajectories, the red and blue curves represent the exponential references with  exponents $-\eta M/2$ and $-\eta M$ respectively. The figures at the bottom are the semi-log versions of the ones at the top.}
      \label{FB3_0_VD_fig}
\end{figure}

\begin{figure}[thpb]
\centering
      \includegraphics[width=12.5cm]{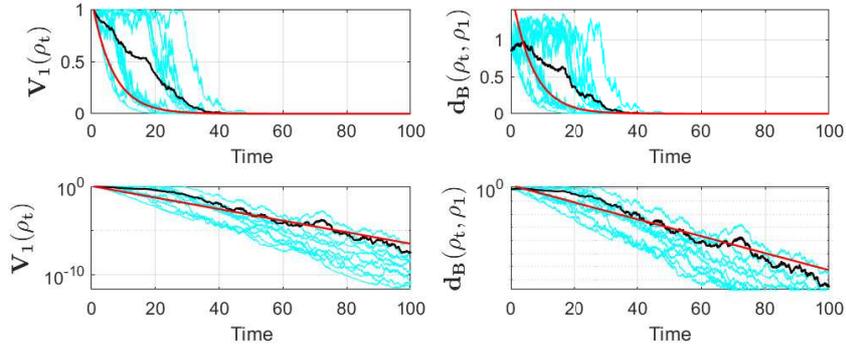}
      \caption{Exponential stabilization of a three-level quantum angular momentum system towards $\boldsymbol \rho_1$ with the feedback law~\eqref{u_t N} starting at $diag(0.3,0.4,0.3)$ with $\omega=0$, $\eta=0.3$, $M=1$, $\alpha = 0.3$, $\beta = 10$: the black curve represents the mean value of 10 arbitrary sample trajectories, the red curve represents the exponential reference with exponent $-\eta M/2$. The figures at the bottom are the semi-log versions of the ones at the top.}
      \label{FB3_1_VD_fig_interior}
\end{figure}

The simulations in the case $u\equiv 0$ are shown in  Fig.~\ref{QSR3fig}. In particular, we observe that the expectation of the Lyapunov function $\mathbb{E}(V(\rho_t))$ is bounded by the exponential function $V(\rho_0)e^{-\frac{\eta M}{2}t}$, and the expectation of the Bures distance $\mathbb{E}(d_B(\rho_t,\bar{E}))$ is always below the exponential function $C_2/C_1\,d_B(\rho_0,\bar{E})e^{-\frac{\eta M}{2}t}$, with $C_1 = 1/2$ and $C_2 = 3$ (see Equation~\eqref{C1d<=V<=C2d}) in accordance with the results of Section~\ref{sec-qsr}.
Next, we set $\boldsymbol \rho_0$ as the target eigenstate; the corresponding simulations  with a feedback law of the form~\eqref{u_t 0 2J}  and initial condition $\boldsymbol \rho_2$ are shown in Fig.~\ref{FB3_0_VD_fig}. For this case, we note that a larger $\alpha$ can speed up the exit of the trajectories  from a neighborhood of the eigenstate $\boldsymbol \rho_2.$ Similarly, a larger $\gamma$ may speed up the accessibility of a neighborhood of the target state $\boldsymbol \rho_0.$ Finally, a larger $\beta$  can  weaken the role of the first term in  the feedback law~\eqref{u_t 0 2J} on neighborhoods of the target state  (a more detailed discussion for the two-level case may be found in~\cite{liang2018exponential}). 
Then, we set $\boldsymbol \rho_1$ as the target eigenstate; the simulations with a feedback law of  the form~\eqref{u_t N} and initial condition $diag(0.3,0.4,0.3)$ (in the interior of $\mathcal S$) are shown in Fig.~\ref{FB3_1_VD_fig_interior}. 
Finally, we repeat the last simulations for the case where the initial condition is  $\boldsymbol \rho_2$. As simulations show, the trajectories enter immediately in the interior of $\mathcal S$ and converge exponentially towards the target state. 
\begin{figure}[thpb]
\centering
      \includegraphics[width=12.5cm]{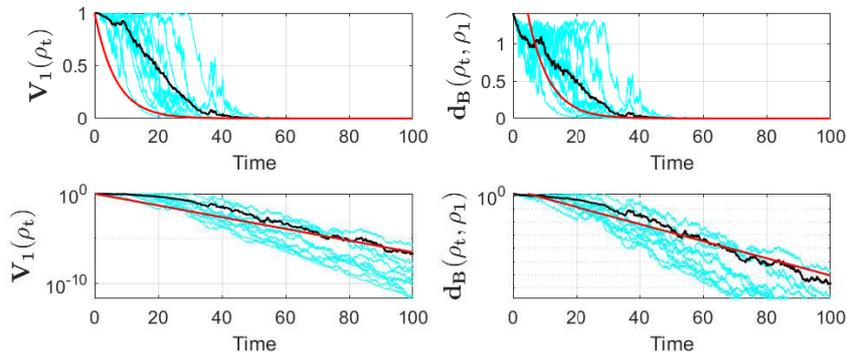}
      \caption{Exponential stabilization of a three-level quantum angular momentum system towards $\boldsymbol \rho_1$ with the feedback law~\eqref{u_t N}starting at $\boldsymbol \rho_2$ with $\omega=0$, $\eta=0.3$, $M=1$, $\alpha = 0.3$, $\beta = 10$: the black curve represents the mean value of  10 arbitrary sample trajectories, the red curve represents the exponential reference with exponent $-\eta M/2$. The figures at the bottom are the semi-log versions of the ones at the top.}
      \label{FB3_1_VD_fig}
\end{figure}

\section{Conclusion and perspectives}
In this paper, we have studied the asymptotic behavior of trajectories associated with quantum angular momentum systems for the cases with and without feedback law. Firstly, for the system with zero control, we have shown the exponential convergence towards the set of eigenstates of the measurement operator $J_z$ (quantum state reduction with exponential rate $\eta M/2$). We next proved  the exponential convergence of $N$-level quantum angular momentum systems towards an arbitrary predetermined target eigenstate under some general conditions on the feedback law. This was obtained by applying stochastic Lyapunov techniques and analyzing the asymptotic behavior of quantum trajectories. For illustration, we have provided a parametrized feedback law satisfying our general conditions to stabilize the system exponentially towards the  target state. 

Further research lines will address  the possibility of extending our results in presence of delays, or for exponential stabilization of entangled states with applications in quantum computing. In particular, alternative choices of the measurement operator may be investigated to prepare predetermined entangled target states, such as Dicke or GHZ states. 

\bibliographystyle{siamplain}
\bibliography{LIANG_SIAM}

\begin{thebibliography}{10}

\bibitem{abe2008analysis}
{\sc T.~Abe, T.~Sasaki, S.~Hara, and K.~Tsumura}, {\em Analysis on behaviors of
  controlled quantum systems via quantum entropy}, IFAC Proceedings Volumes, 41
  (2008), pp.~3695--3700.

\bibitem{adler2001martingale}
{\sc S.~L. Adler, D.~C. Brody, T.~A. Brun, and L.~P. Hughston}, {\em Martingale
  models for quantum state reduction}, Journal of Physics A: Mathematical and
  General, 34 (2001), p.~8795.

\bibitem{ahn2002continuous}
{\sc C.~Ahn, A.~C. Doherty, and A.~J. Landahl}, {\em Continuous quantum error
  correction via quantum feedback control}, Physical Review A, 65 (2002),
  p.~042301.

\bibitem{amini2011design}
{\sc H.~Amini, P.~Rouchon, and M.~Mirrahimi}, {\em Design of strict
  control-lyapunov functions for quantum systems with {QND} measurements}, in
  50th IEEE Conference on Decision and Control and European Control Conference
  (CDC-ECC), 2011, 2011, pp.~8193--8198.

\bibitem{amini2013feedback}
{\sc H.~Amini, R.~A. Somaraju, I.~Dotsenko, C.~Sayrin, M.~Mirrahimi, and
  P.~Rouchon}, {\em Feedback stabilization of discrete-time quantum systems
  subject to non-demolition measurements with imperfections and delays},
  Automatica, 49 (2013), pp.~2683--2692.

\bibitem{armen2002adaptive}
{\sc M.~A. Armen, J.~K. Au, J.~K. Stockton, A.~C. Doherty, and H.~Mabuchi},
  {\em Adaptive homodyne measurement of optical phase}, Physical Review
  Letters, 89 (2002), p.~133602.

\bibitem{belavkin1983theory}
{\sc V.~P. Belavkin}, {\em On the theory of controlling observable quantum
  systems}, Avtomatika i Telemekhanika,  (1983), pp.~50--63.

\bibitem{belavkin1989nondemolition}
{\sc V.~P. Belavkin}, {\em Nondemolition measurements, nonlinear filtering and
  dynamic programming of quantum stochastic processes}, in Modeling and Control
  of Systems, Springer, 1989, pp.~245--265.

\bibitem{belavkin1992quantum}
{\sc V.~P. Belavkin}, {\em Quantum stochastic calculus and quantum nonlinear
  filtering}, Journal of Multivariate analysis, 42 (1992), pp.~171--201.

\bibitem{belavkin1995quantum}
{\sc V.~P. Belavkin}, {\em Quantum filtering of markov signals with white
  quantum noise}, in Quantum communications and measurement, Springer, 1995,
  pp.~381--391.

\bibitem{bengtsson2017geometry}
{\sc I.~Bengtsson and K.~{\.Z}yczkowski}, {\em Geometry of quantum states: an
  introduction to quantum entanglement}, Cambridge University Press, 2017.

\bibitem{bouten2008separation}
{\sc L.~Bouten and R.~Van~Handel}, {\em On the separation principle in quantum
  control}, in Quantum stochastics and information: statistics, filtering and
  control, World Scientific, 2008, pp.~206--238.

\bibitem{bouten2007introduction}
{\sc L.~Bouten, R.~Van~Handel, and M.~R. James}, {\em An introduction to
  quantum filtering}, SIAM Journal on Control and Optimization, 46 (2007),
  pp.~2199--2241.

\bibitem{cardona2018exponential}
{\sc G.~Cardona, A.~Sarlette, and P.~Rouchon}, {\em Exponential stochastic
  stabilization of a two-level quantum system via strict lyapunov control}, in
  IEEE Conference on Decision and Control, 2018, pp.~6591--6596.

\bibitem{davies1969quantum}
{\sc E.~B. Davies}, {\em Quantum stochastic processes}, Communications in
  Mathematical Physics, 15 (1969), pp.~277--304.

\bibitem{davies1976quantum}
{\sc E.~B. Davies}, {\em Quantum theory of open systems}, Academic Press, 1976.

\bibitem{dotsenko2009quantum}
{\sc I.~Dotsenko, M.~Mirrahimi, M.~Brune, S.~Haroche, J.-M. Raimond, and
  P.~Rouchon}, {\em Quantum feedback by discrete quantum nondemolition
  measurements: Towards on-demand generation of photon-number states}, Physical
  Review A, 80 (2009), p.~013805.

\bibitem{dynkin1965markov}
{\sc E.~B. Dynkin}, {\em Markov processes}, in Markov Processes, Springer,
  1965, pp.~77--104.

\bibitem{hudson1984quantum}
{\sc R.~L. Hudson and K.~R. Parthasarathy}, {\em Quantum {I}to's formula and
  stochastic evolutions}, Communications in Mathematical Physics, 93 (1984),
  pp.~301--323.

\bibitem{kato1976perturbation}
{\sc T.~Kato}, {\em Perturbation theory for linear operators}, vol.~132,
  Springer, 1976.

\bibitem{khasminskii2011stochastic}
{\sc R.~Khasminskii}, {\em Stochastic stability of differential equations},
  vol.~66, Springer, 2011.

\bibitem{liang2018exponential}
{\sc W.~Liang, N.~H. Amini, and P.~Mason}, {\em On exponential stabilization of
  spin-$\frac12$ systems}, in IEEE Conference on Decision and Control, 2018,
  pp.~6602--6607.

\bibitem{mabuchi2005principles}
{\sc H.~Mabuchi and N.~Khaneja}, {\em Principles and applications of control in
  quantum systems}, International Journal of Robust and Nonlinear Control:
  IFAC-Affiliated Journal, 15 (2005), pp.~647--667.

\bibitem{mao2007stochastic}
{\sc X.~Mao}, {\em Stochastic differential equations and applications},
  Woodhead Publishing, 2007.

\bibitem{mirrahimi2009feedback}
{\sc M.~Mirrahimi, I.~Dotsenko, and P.~Rouchon}, {\em Feedback generation of
  quantum fock states by discrete qnd measures}, in IEEE Conference on Decision
  and Control, 2009, pp.~1451--1456.

\bibitem{mirrahimi2007stabilizing}
{\sc M.~Mirrahimi and R.~Van~Handel}, {\em Stabilizing feedback controls for
  quantum systems}, SIAM Journal on Control and Optimization, 46 (2007),
  pp.~445--467.

\bibitem{pellegrini2008existence}
{\sc C.~Pellegrini}, {\em Existence, uniqueness and approximation of a
  stochastic schr{\"o}dinger equation: the diffusive case}, The Annals of
  Probability,  (2008), pp.~2332--2353.

\bibitem{protter2004stochastic}
{\sc P.~E. Protter}, {\em Stochastic integration and differential equations},
  2004.

\bibitem{revuz2013continuous}
{\sc D.~Revuz and M.~Yor}, {\em Continuous martingales and Brownian motion},
  vol.~293, Springer, 2013.

\bibitem{rogers2000diffusions1}
{\sc L.~G. Rogers and D.~Williams}, {\em Diffusions, Markov processes and
  martingales: Volume 1, Foundations}, vol.~1, Cambridge University Press,
  2000.

\bibitem{rogers2000diffusions2}
{\sc L.~G. Rogers and D.~Williams}, {\em Diffusions, Markov processes and
  martingales: Volume 2, It{\^o} calculus}, vol.~2, Cambridge university press,
  2000.

\bibitem{sarlette2017deterministic}
{\sc A.~Sarlette and P.~Rouchon}, {\em Deterministic submanifolds and analytic
  solution of the quantum stochastic differential master equation describing a
  monitored qubit}, Journal of Mathematical Physics, 58 (2017), p.~062106.

\bibitem{sayrin2011real}
{\sc C.~Sayrin, I.~Dotsenko, X.~Zhou, B.~Peaudecerf, T.~Rybarczyk, S.~Gleyzes,
  P.~Rouchon, M.~Mirrahimi, H.~Amini, M.~Brune, J.-M. Raimond, and S.~Haroche},
  {\em Real-time quantum feedback prepares and stabilizes photon number
  states}, Nature, 477 (2011), pp.~73--77.

\bibitem{stroock1972support}
{\sc D.~W. Stroock and S.~R. Varadhan}, {\em On the support of diffusion
  processes with applications to the strong maximum principle}, in Proceedings
  of the Sixth Berkeley Symposium on Mathematical Statistics and Probability
  (Univ. California, Berkeley, Calif., 1970/1971), vol.~3, 1972, pp.~333--359.

\bibitem{tsumura2008global}
{\sc K.~Tsumura}, {\em Global stabilization at arbitrary eigenstates of
  n-dimensional quantum spin systems via continuous feedback}, in American
  Control Conference, 2008, 2008, pp.~4148--4153.

\bibitem{van2005feedback}
{\sc R.~Van~Handel, J.~K. Stockton, and H.~Mabuchi}, {\em Feedback control of
  quantum state reduction}, IEEE Transactions on Automatic Control, 50 (2005),
  pp.~768--780.

\bibitem{xiong2008introduction}
{\sc J.~Xiong}, {\em An introduction to stochastic filtering theory}, vol.~18,
  Oxford University Press, 2008.

\bibitem{yamamoto2007feedback}
{\sc N.~Yamamoto, K.~Tsumura, and S.~Hara}, {\em Feedback control of quantum
  entanglement in a two-spin system}, Automatica, 43 (2007), pp.~981--992.

\end{thebibliography}
\end{document}